\documentclass[final]{siamart1116}
\usepackage{graphicx,helvet}
\usepackage{epstopdf}
\usepackage[percent]{overpic}
\usepackage{xcolor}
\usepackage{multirow}

\usepackage{epsfig}
\usepackage{cite}
\usepackage[english]{babel}
\usepackage[utf8]{inputenc}
\usepackage{amssymb}

\allowdisplaybreaks

\newcommand{\eps}{\varepsilon}
\renewcommand{\epsilon}{\varepsilon}

\newcommand{\R}{\mathbb{R}}

\newcounter{dctr}[section]
\numberwithin{equation}{section}

\newtheorem{deff}{Definition}[section]
\newtheorem{remark}{Remark}[section]

\newcommand{\bigO}{\mathcal{O}}

\newcommand{\C}{\mathcal{C}}

\newcommand{\Z}{\mathbb{Z}}
\newcommand{\N}{\mathbb{N}}

\def\sf#1#2{\sum_{i} \frac{#1_i #2_i}{v_i-\lambda}}

\headers{Third-order WENO scheme with optimal accuracy}{Baeza, B\"{u}rger, Mulet, and Zor\'{\i}o}

\allowdisplaybreaks
\numberwithin{equation}{section}

\title{An efficient third-order WENO scheme with unconditionally optimal accuracy} 

\date{\today}

\author{Antonio Baeza\thanks{$^{\mathrm{A}}$Departament de Matem\`{a}tiques, 
Universitat de Val\`{e}ncia, Av.\ Vicent Andr\'es Estell\'es, E-46100
Burjassot,  
 Spain.  E-Mail: 
   {\tt antonio.baeza@uv.es}}
   \and 
 Raimund B\"{u}rger\thanks{$^{\mathrm{B}}$CI$^{\mathrm{2}}$MA and 
Departamento de Ingenier\'{\i}a Matem\'{a}tica,   Universidad de
Concepci\'{o}n, Casilla \mbox{160-C}, Concepci\'{o}n, Chile.  E-Mail:
  {\tt rburger@ing-mat.udec.cl}}
  \and 
   Pep Mulet\thanks{$^{\mathrm{C}}$Departament de Matem\`{a}tiques, 
Universitat de Val\`{e}ncia, Av.\ Vicent Andr\'es Estell\'es, E-46100
Burjassot,  
 Spain.  E-Mail: 
   {\tt pep.mulet@uv.es}}
   \and 
  David Zor\'{\i}o\thanks{$^{\mathrm{D}}$CI$^{\mathrm{2}}$MA,   Universidad de
Concepci\'{o}n, Casilla \mbox{160-C}, Concepci\'{o}n, Chile.  E-Mail:
  {\tt dzorio@ci2ma.udec.cl}}}

\usepackage[normalem]{ulem}
\begin{document}

\maketitle 

\begin{abstract}
  A novel scheme, based on third-order Weighted Essentially Non-Oscillatory (WENO) reconstructions, is presented. It attains unconditionally optimal accuracy when the data is smooth enough, even in presence of critical points, and second-order accuracy if a discontinuity crosses the data. The key to attribute these properties to this scheme is the inclusion of an additional node in the data stencil, which is only used in the computation of the weights measuring the smoothness. The accuracy properties of this scheme are proven in detail and several numerical experiments are presented, which show that this scheme is more efficient in terms of the error reduction versus CPU time than its traditional third-order counterparts as well as several higher-order WENO schemes that are found in the literature.
\end{abstract} 

\date{\today}



\section{Introduction}

\subsection{Scope} \label{subsec:scope}

Weighted Essentially Non-Oscillatory (WENO) schemes have become very
popular, especially in the context of hyperbolic conservation laws,
since they were proposed in \cite{LiuOsherChan94} and later improved
in \cite{JiangShu96}. One of the most used schemes in the literature
is the fifth-order WENO scheme, which in general
attains satisfactory results on weak solutions of hyperbolic conservation laws.

Albeit traditional third order methods are also widely used,
  the accuracy loss near smooth extrema is an issue that lowers
  significantly the accuracy of the numerical solution, even for
  problems with weak solutions, in which it is significantly smeared.

In this paper we inspect the causes of the misperformance involving the traditional third-order WENO schemes through an analysis of their accuracy near critical points. We propose several solutions to this issue, by first proving that it is impossible to prevent accuracy loss near critical points in stencils with only three points, and then showing that it is possible to do so with stencils of at least four points. Ultimately the goal is to present a genuine third-order scheme that is competitive with the most widely used fifth-order schemes for problems with weak solutions.

\subsection{Related work} To put this work into the proper
perspective, we mention several previous attempts that have been made
in order to solve the issue involving the accuracy loss near critical
points. For instance, in \cite{WuLiangZhao}, the authors propose a
novel smoothness measure based on introducing an additional exponent
in the weight formula proposed in \cite{WuZhao} associated to WENO-N3
schemes. However, although this measure  solves  the issue of the accuracy loss near
critical points, the resulting weights 
depend on the scaling of the data due to the additional
exponent. Other works improving this idea have been also done, but
the issue of the weights depending on the scaling of the data
in the weight design still remains; see for instance
\cite{GandeRathodRathan17,GandeRathodRathan18,XuWu}.

Many other works deal with the issue by tuning the
parameter $\varepsilon$ appearing in the weight design, which was
initially conceived to be a small quantity used to avoid divisions by
zero near constant data, but that was proven later to be crucial to avoid
the accuracy loss near critical points if it was scaled
properly, see for instance \cite{SINUM2011}. Some recent
works which deal with this issue in third-order schemes and limiters, for
instance \cite{Schmidtmann16,SchmidtmannJSC}.

In the case of higher accuracy order methods, the issue
    of the accuracy loss near critical points, without relying on
    tuning or scaling parameters, has been handled more
    broadly in the literature by proposing new weight designs, such as
 the WENO-M \cite{HenrickAslamPowers2005}, WENO-Z \cite{BorgesCarmonaCostaEtAl2008} 
    and Yamaleev-Carpenter methods \cite{YamaleevCarpenter2009b},
   obtaining partial solutions to the problem for schemes of arbitrary order.
   In \cite{BBMZSINUM} we proposed a method that completely solves the issue 
     for schemes of order higher than 3. In fact, the present work can be seen
  as a complement of \cite{BBMZSINUM}, in which
  the third-order case, that cannot be fit in the general framework, is separately tackled 
  through a new approach.

We will show in this work that it is not possible to build a third-order
reconstruction with a stencil of three points satisfying at once the
following properties: 
\begin{itemize}
\item Detection of discontinuities in the data.
\item Detection of critical points in the data. 
\item Independence of non-linear weights  of the scaling of the data (the
  issue appearing in \cite{WuLiangZhao}).
\item  Unnecessity of tuning/scaling~$\varepsilon$ (in contrast to the proposals  
  in, e.g.,  \cite{SINUM2011,SchmidtmannJSC}).
\end{itemize}
And, once exposed, we will propose a novel WENO3 reconstruction method
satisfying at once the aforementioned properties, by using stencils
with an additional point, namely, a stencil containing a total of four
points. This does not represent an increase of the stencil used to
compute the numerical divergence in a semi-discrete scheme, as proven
in Section \ref{subsec:dependence}.

\subsection{Outline of the paper}

This paper is organized as follows: Section \ref{sec:oweno3} starts with
some preliminaries and definitions that will be used along the work,
presented in Subsections \ref{subsec:prel} and \ref{subsec:wenorec},
followed by a motivation in Subsection \ref{subsec:weno3} in which we
prove through a counterexample that a third-order WENO scheme cannot
attain the optimal accuracy near critical points if a stencil of only
three points is used, but that it is possible to attain the optimal
accuracy even near critical points if an additional point is
added. The proposed scheme, attaining unconditionally third order, is
presented in Subsection \ref{subsec:oweno3}. In Section \ref{sec:wenohyp} the
key to use this reconstruction strategy in the context of third-order
schemes for hyperbolic conservation laws without increasing the
computational domain is shown. Section \ref{sec:num_exp}
stands for several validation numerical experiments in which our
proposed schemes are compared against the most commonly used
fifth-order scheme in terms of efficiency; finally, in Section
\ref{sec:conclusions} some conclusions are drawn.

\section{Optimal third-order scheme}\label{sec:oweno3}

\subsection{Preliminaries}\label{subsec:prel}
\begin{deff}
  Assume that $\alpha\in\mathbb{Z}$.  
   We write 
  $f(h) =\bigO(h^\alpha)$ to denote that  $\limsup_{h\to 0} |f(h)/ h^\alpha | <\infty$, 
   and $f(h) =\bar{\bigO}(h^\alpha)$ if $f(h) =\bigO(h^\alpha)$ and in addition $\liminf_{h\to 0}| f(h) / h^\alpha | > 0$.  
\end{deff} 
Since, for positive functions $f$ and $g$, 
\begin{align*}
  \limsup_{h\to 0} f(h)g(h)&\leq \limsup_{h\to 0} f(h) \limsup_{h\to
    0}  g(h),\\
 \liminf_{h\to 0} f(h)g(h)&\geq \liminf_{h\to 0} f(h) \liminf_{h\to 0} g(h),
\end{align*}
it follows that
$\bigO(h^{\alpha}) \bigO(h^{\beta}) = \bigO(h^{\alpha+\beta})$ and 
$\smash{\bar\bigO(h^{\alpha}) \bar\bigO(h^{\beta}) =
\bar\bigO(h^{\alpha+\beta})}$.

\subsection{Third-order WENO reconstructions}\label{subsec:wenorec}

For the sake of exposition we briefly describe two classical 
third-order WENO approaches. The first  is  the third-order WENO method
defined by   the Jiang-Shu approach \cite{JiangShu96} (henceforth,
JS-WENO3) and the second is 
the third-order WENO method through the Yamaleev-Carpenter approach
\cite{YamaleevCarpenter2009a,YamaleevCarpenter2009b} (henceforth,
YC-WENO3). Since they have
many  parts in common, we will describe both approaches altogether while
pointing out the key differences when necessary.

The input for both cases is an equally-spaced three-point stencil
$(x_{-1},x_0,x_1)$, $x_i-x_{i-1}=h>0$, $i\in\{0,1\}$, associated with 
values $(f_{-1},f_{0},f_{1})$, where either 
  $f_i=f(x_i)$ (reconstructions from point
values, namely, an interpolation procedure
  in which the data from the stencil is interpreted as point values
  of a function, with the reconstruction being a point value of that
  function) or  
\begin{align*} 
f_i=\frac{1}{h}\int_{x_{i-1/2}}^{x_{i+1/2}}f(x)\, \mathrm{d}x
\end{align*} 
(reconstructions from cell averages, namely, an interpolation
  procedure
  in which the data from the stencil is interpreted as cell averages
  of a function, with the reconstruction being a point value of that
  function). Here  
$x_{i+1/2}= (x_i+x_{i+1})/2$, and $\varepsilon>0$ is    a parameter 
whose original purpose is to be merely a small positive quantity avoiding
divisions by zero. We assume that a right-biased reconstruction is sought, so the output
is intended to be an approximation of $f(x_{1/2})$.
 The smoothness indicators \cite{JiangShu96} are then defined as follows:
\begin{align}  \label{eq:I0I1} 
  I_0 := (f_0-f_{-1})^2,\quad 
  I_1 := (f_1-f_0)^2,
\end{align} 
along with the corresponding interpolating polynomials associated to
each 2-point substencil:
\begin{align}   \label{p0xhalf} 
  p_0(x_{1/2}) = -\frac{1}{2} f_{-1} + \frac{3}{2} f_0,\quad 
  p_1(x_{1/2}) = \frac{1}{2} f_0 + \frac{1}{2} f_1.
\end{align} 
Now, in each case we define 
\begin{equation}\label{eq:alphai}\alpha_i=
\begin{cases}
 \displaystyle  \frac{c_i}{I_i+\varepsilon} & \textnormal{for the
   JS-WENO3 method,}\\[4mm]
  c_i\left( \displaystyle 1+\frac{\sigma}{I_i+\varepsilon}\right) &
  \textnormal{for the YC-WENO3 method,}
\end{cases}
\end{equation}
where $\sigma=(f_1-2f_0+f_{-1})^2$ and
\begin{align}  \label{c0c1def} 
(c_0, c_1) = \begin{cases} 
 (1/4, 3/4) & \text{in case of reconstructions from point values,} \\
 (1/3, 2/3) & \text{in case of reconstructions from cell averages.} 
 \end{cases}
\end{align} 
Then, the non-linear weights are computed as
$$\omega_i=\frac{\alpha_i}{\alpha_0+\alpha_1}, \quad i=0,1,$$
and the WENO reconstruction is finally given by
$$p_2(x_{1/2})= \omega_0p_0(x_{1/2})+\omega_1p_1(x_{1/2}).$$

\begin{remark}
Although the denominator appearing in the expressions of $\alpha_i$ in \eqref{eq:alphai} is commonly chosen as $(I_i+\varepsilon)^2$ for third-order schemes, in \cite[Note 2]{SINUM2011} it was proven that, in general, for a ($2r-1$)-th order scheme, a sufficient condition to attain the suboptimal $r$-th order accuracy when a discontinuity crosses the stencil is using a denominator of the form $(I_i+\varepsilon)^p$ with $2p\geq r$. Therefore, in the particular case of third order ($r=2$), it suffices to choose $p=1$, in which theoretical order properties identical to the case $p=2$ are attained.
\end{remark}

\subsection{On the accuracy loss of third-order WENO schemes}\label{subsec:weno3}

In the subsequent text we will abuse language by referring to the values
of a function on a stencil as the stencil itself.

The following example shows that if a grid $x_{i,h}=z+(c+i)h$, $i\in\{-1,0,1\}$, samples a function $f\in\C^{2}$ such that $f'(z)=0$ and $f''(z)\neq0$, then there are cases in which, ignoring the scaling of the stencil 
 $(f_{-1,h},f_{0,h},f_{1,h})$, with $f_{i,h}=f(x_{i,h})$, the
 reconstruction obtained from that stencil for any given $h$ is the
 same as if the function had a discontinuity in it, and thus there
 cannot be scaling-independent and dimensionless parameters
 constructed from the data capable of distinguishing one case from the
 other.

Let us consider, on one hand, an extreme case by considering $f:\R\to\R$ given by $f(x)=4x^2$, which satisfies $f'(0)=0$ and $f''(0)=8\neq0$, and the grid $x_{i,h}=(\frac{1}{2}+i)h$, $i\in\{-1,0,1\}$. Then the stencil $F_h=(f_{-1,h},f_{0,h},f_{1,h})$, with $f_{i,h}=f(x_{i,h})$, is given by $F_h=(h^2,h^2,9h^2)$.

On the other hand, we define $g:\R\to\R$ given by $g(x)=1$ if $x<0$ and $g(x)=9$ if $x\geq0$, with the same grid as above. Then the stencil $G_h=(g_{-1,h},g_{0,h},g_{1,h} )$ is given by $G_h=(1,1,9)$.

Now, the relationship $F_h=h^2G_h$ holds for all $h>0$; that is, both stencils are, for fixed $h$, a scaled version of the same stencil. Therefore, any procedure to analyze smoothness agnostic about the scaling of the data will fail at distinguishing  the first case, consisting of smooth data, from the second one, based on data taken from both sides of a discontinuity. Therefore, such procedure, depending on its construction, will either detect asymptotically both cases as smooth data, or will interpret both as discontinuous data, being in both cases wrong (giving either false negatives or false positives).  The traditional third-order WENO schemes,  
 belong to  the latter group, in which the detection of discontinuities is prioritized against the detection of critical points, and thus the latter ones are interpreted incorrectly as discontinuities.

\subsection{Unconditionally optimal third-order scheme with an additional node}\label{subsec:oweno3}

We next present a novel scheme with essentially non-oscillatory
properties which attains unconditionally the optimal order of
accuracy.

Let $S:= (f_{-1},f_0,f_1 )$ be a stencil from a uniform grid, $f_i=f(x_i)$, $i\in\{-1,0,1\}$, $x_i=z+(c+i)h$, $i\in\{-1,0,1\}$, and $\bar{S}:=S\cup\{f_2\}= (f_{-1},f_0,f_1,f_2 )$ the extended stencil. Let us assume that one wishes to perform a (right-biased with respect to $S$) reconstruction at $x_{1/2}:= (x_0+x_1)/2 $ accounting for discontinuities. Then, both for reconstructions from point values and from cell averages, we define the following items:

We define the corresponding interpolating polynomials associated to the substencils $S_0= (f_{-1},f_0 )$ and $S_1= (f_0,f_1 )$ evaluated at $x_{1/2}$, which are given by \eqref{p0xhalf}.  
Their associated Jiang-Shu smoothness indicators are thus given by 
 \eqref{eq:I0I1}. 
One of the keys here is to define also an additional smoothness
indicator, in which the additional node is used, namely
\begin{align} \label{eq:I2} 
  I_2 := (f_2-f_1)^2.
\end{align} 
Now, given a small quantity $\varepsilon>0$, we define the weights
\begin{align} \label{eq:tw}  
  \tilde{\omega}_0 := \frac{I_1+\varepsilon}{I_0+I_1+2\varepsilon},\quad 
  \tilde{\omega}_1 := \frac{I_0+\varepsilon}{I_0+I_1+2\varepsilon}=1-\tilde{\omega}_0.
\end{align} 
We introduce now the corrector weight, given by
\begin{align} \label{eq:cw} 
\omega=\frac{J}{J+\tau+\varepsilon} \quad 
\text{with $J=I_0(I_1+I_2)+(I_0+I_1)I_2$,} 
\end{align} 
which clearly satisfies $0\leq\omega\leq1$, and $\tau$ the product of the square of the undivided difference associated to the extended stencil $\bar{S}$ with the sum of the smoothness indicators:
\begin{align} \label{eq:taudef} 
  \tau:=dI,\quad d:=(-f_{-1}+3f_0-3f_1+f_2)^2,\quad I:=I_0+I_1+I_2. 
\end{align} 
We then define the corrected weights as
\begin{align} \label{correctedweights}  
&   \omega_0 := \omega c_0+(1-\omega)\tilde{\omega}_0,\quad 
  \omega_1 :=  \omega c_1+(1-\omega)\tilde{\omega}_1,
 \end{align} 
 where $c_0$ and~$c_1$ are specified in    \eqref{c0c1def}.  
Finally, the reconstruction result is given by
\begin{align} \label{pxhalf} 
 p(x_{1/2} )=\omega_0p_0(x_{1/2})+\omega_1p_1(x_{1/2}).
 \end{align} 

The key to analyze the accuracy of our proposed scheme is to first
study the accuracy of the corrector
weight~$\omega$.

\begin{definition}
  We say that a function $f$ has a  critical point of 
order~$k\geq0$ at $x$ if $\smash{f^{(l)}}(x)=0$ for 
$l=1,\dots,k$ and $\smash{f^{(k+1)}}(x)\neq 0$.
\end{definition}

\begin{proposition}\label{w_accuracy}
  If $f$ has a critical point at $z$ of order $k$, $k\in\{0,1\}$, there holds
  \[
  \omega=\begin{cases}
  1+\bigO(h^{4-2k})+\bigO(\varepsilon) & \textnormal{if }\bar{S}\textnormal{ is smooth, }f\in\C^3,\\
  \bigO(h^2)+\bigO(\varepsilon) & \textnormal{if a discontinuity crosses }S.
  \end{cases}
  \]
\end{proposition}

\begin{proof}
  Clearly, by definition and the fact that $J,\tau\geq0$, there holds $0\leq\omega\leq 1$.
  
  Let us first assume that $\bar{S}$ is smooth with
  $k\in\{0,1\}$. Then, according to \cite[Lemma 2]{BaezaBurgerMuletZorio2018}, if $k=0$, $I_{2,i}=\bar{\mathcal{O}}(h^2)$, $i\in\{0,1,2\}$, and if $k=1$, then 
   there exists $i_0 \in \{ 0,1,2 \}$  such that $I_{2,i_0}=\bar{\mathcal{O}}(h^s)$, for some $s\in\{4,5,6,\ldots\}$, and $I_{2,i}=\bar{\mathcal{O}}(h^4)$, for $i\in\{0,1,2\}$, $i\neq i_0$.

  Therefore, combining these properties, we deduce that $I_{2,0}+I_{2,1}=\bar{\mathcal{O}}(h^{2+2k})$ and that $I_{2,1}+I_{2,2}=\bar{\mathcal{O}}(h^{2+2k})$. Moreover, since either $I_{2,0}=\bar{\mathcal{O}}(h^{2+2k})$ or $I_{2,2}=\bar{\mathcal{O}}(h^{2+2k})$, it can be concluded that
  $$J=I_{2,0}(I_{2,1}+I_{2,2})+(I_{2,0}+I_{2,1})I_{2,2}=\bar{\mathcal{O}}(h^{4+4k})+\bigO(\varepsilon).$$
  On the other hand, 
    \begin{align*}
  d&=(-f_{-1}+3f_0-3f_1+f_2)^2=\bigO(h^6)=\bigO(h^6)+\bigO(\varepsilon),\\
  I&=I_0+I_1+I_2=\bigO(h^{2+2k})=\bigO(h^{2+2k})+\bigO(\varepsilon).
  \end{align*}
  Therefore $\tau=\bigO(h^{8+2k})$. Hence, and since by assumption $J\neq0$,
  \begin{align*}
    \omega&=\frac{J}{J+\tau+\varepsilon}=\frac{1}{\displaystyle1+\frac{\tau}{J}}-\bigO(\varepsilon)=\frac{1}{\displaystyle1+\frac{\bigO(h^{8+2k})}{\bar{\mathcal{O}}(h^{4+4k})}}-\bigO(\varepsilon)=\frac{1}{1+\bigO(h^{4-2k})}-\bigO(\varepsilon)\\
    &=1-\bigO(h^{4-2k})-\bigO(\varepsilon).
  \end{align*}
  
  Finally, let us assume that a discontinuity crosses $S$. Then there exists  $i_0\in\{0,1\}$ such that $I_{2,i_0}=\bar{\mathcal{O}}(1)$. On the other hand, $I_{2,|1-i_0|}=\bar{\mathcal{O}}(h^{2m_0})$ and $I_{2,2}=\bar{\mathcal{O}}(h^{2m_1})$ for some $1\leq m_0,m_1\in\{1,2,3,\dots\}$. Now, by these considerations, we have
  \begin{align*}
    I_{2,0}(I_{2,1}+I_{2,2})&=\begin{cases}
    \bar\bigO(h^m) & \text{if $i_0=0$,}  \\
    \bar\bigO(h^{m_0}) & \text{if $i_0=1$,} 
    \end{cases}\quad 
    (I_{2,0}+I_{2,1})I_{2,2} =\bar\bigO(h^{m_1})  
  \end{align*}
  with $m:=\max\{m_0,m_1\}$.

  Under any of these combinations, we obtain
  $$J=I_{2,0}(I_{2,1}+I_{2,2})+(I_{2,0}+I_{2,1})I_{2,2}=\bar{\mathcal{O}}(h^{2m}).$$

  On the other hand, since in this case there holds
  \begin{align*}
  d&=(-f_{-1}+3f_0-3f_1+f_2)^2=\bar{\mathcal{O}}(1)=\bar{\mathcal{O}}(1)+\bigO(\varepsilon),\\
  I&=I_0+I_1+I_2=\bar{\mathcal{O}}(1)=\bar{\mathcal{O}}(1)+\bigO(\varepsilon),
  \end{align*}
  then $\tau=\bar{\mathcal{O}}(1)$ and
  \begin{align*}
    \omega&=\frac{J}{J+\tau+\varepsilon}=\frac{1}{\displaystyle1+\frac{\tau}{J}}-\bigO(\varepsilon)=\frac{1}{\displaystyle1+\frac{\bar{\mathcal{O}}(1)^2}{\bar{\mathcal{O}}(h^{2m})}}-\bigO(\varepsilon)=\frac{1}{\displaystyle1+\frac{\bar{\mathcal{O}}(1)}{\bar{\mathcal{O}}(h^{2m})}}-\bigO(\varepsilon)\\
    &=\frac{1}{\displaystyle1+\bar{\mathcal{O}}(h^{-2m})}-\bigO(\varepsilon)=\frac{1}{\displaystyle\bar{\mathcal{O}}(h^{-2m})}-\bigO(\varepsilon)=\bar{\mathcal{O}}(h^{2m})+\bigO(\varepsilon)=\bigO(h^2)+\bigO(\varepsilon),
  \end{align*}
  which completes the proof.
\end{proof}

Now, let us focus on the computation of the corrected weights.

\begin{proposition}\label{wi_accuracy}
  For $i\in\{0,1\}$ there holds
  \[\omega_i=
  \begin{cases}
    c_i+\bigO(h^{4-2k})+\bigO(\varepsilon) & \textnormal{if }\bar{S}\textnormal{ contains smooth data,}\\
    \bigO(h^2)+\bigO(\varepsilon) & \textnormal{if a discontinuity crosses }S_i,\\
    \bigO(1)+\bigO(\varepsilon) & \textnormal{if a discontinuity crosses }S\textnormal{, but not }S_i.
  \end{cases}
  \]
\end{proposition}

\begin{proof}
  We first recall that
  $\omega_i=\omega c_i+(1-\omega)\tilde{\omega}_i$. 
  If  $\omega=1-\bigO(h^{m_0})-\bigO(\varepsilon)$ for some $m_0\geq0$, then
  $$\omega_i=\bigl(1-\bigO(h^{m_0})-\bigO(\varepsilon)\bigr)c_i+\bigl(\bigO(h^{m_0})+\bigO(\varepsilon)\bigr)\tilde{\omega}_i=c_i+\bigO(h^{m_0})+\bigO(\varepsilon),$$
  where we have taken into account that $\tilde{\omega}_i$ is an expression at most $\bigO(1)$, since in particular $0\leq\tilde{\omega}_i\leq1$.
  Therefore, using Proposition~\ref{w_accuracy}, we obtain  the result.

  On the other hand, if $\omega$ satisfies $\omega=\bigO(h^{2m_1})+\bigO(\varepsilon)$ for some $m_1\geq1$, then
  $$\omega_i= \bigl(\bigO(h^{2m_1})+\bigO(\varepsilon) \bigr)c_i+ \bigl(1-\bigO(h^{2m_1})-\bigO(\varepsilon) \bigr)\tilde{\omega}_i=\tilde{\omega}_i+\bigO(h^{2m_1})+\bigO(\varepsilon).$$
  Hence, in this case we must focus on the analysis of the accuracy for $\omega_i$. By Proposition~\ref{w_accuracy} we have that $\omega=\bigO(h^{2m_1})$, $m_1>0$, if a discontinuity crosses $S$.

  In such case, there  exists $i_0 \in\{0,1\}$ such that $I_{i_0}=\bar{\mathcal{O}}(1)$, whereas $I_{1-i_0}=\bigO(h^2)$. Therefore, in this case we have
  \begin{align*}
    \tilde{\omega}_{i_0}&=\frac{I_{1-i_0}+\varepsilon}{I_0+I_1+2\varepsilon}=\frac{\bigO(h^2)}{\bar{\mathcal{O}}(1)}+\bigO(\varepsilon)=\bigO(h^2)+\bigO(\varepsilon),\\
    \tilde{\omega}_{1-i_0}&=\frac{I_{i_0}+\varepsilon}{I_0+I_1+2\varepsilon}=\frac{\bar{\mathcal{O}}(1)}{\bar{\mathcal{O}}(1)}+\bigO(\varepsilon)=\bigO(1)+\bigO(\varepsilon).
  \end{align*}
  Therefore, taking into account that $\omega_i=\tilde{\omega}_i+\bigO(h^{2m_1})+\bigO(\varepsilon)$ with $m_1\geq1$, we obtain 
  \begin{align*}
    \omega_{i_0} =\bigO(h^2)+\bigO(\varepsilon),\quad 
    \omega_{1-i_0} =\bigO(1)+\bigO(\varepsilon).
  \end{align*}
\end{proof}

\begin{theorem}\label{accuracy}
 The reconstruction $p(x_{1/2})$ satisfies 
  \[
  p(x_{1/2})=\begin{cases}
  f(x_{1/2})+\bigO(h^3) & \textnormal{if }\bar{S}\textnormal{ is smooth},\\
  f(x_{1/2})+\bigO(h^2) & \textnormal{if a discontinuity crosses }S.
  \end{cases}
  \]
\end{theorem}
\begin{proof}
  This is a direct consequence of the application of Proposition \ref{wi_accuracy} to the expression~\eqref{pxhalf},  
  where we also take into account that both for reconstructions from point
  values and from cell averages, the ideal weights~$c_i$,
  $i\in\{0,1\}$, satisfy that $c_0p_0(x_{1/2})+c_1p_1(x_{1/2})$ equals
   the corresponding third-order reconstruction of the same type at~$x_{1/2}$.
\end{proof}

\begin{remark}
  The cases in which the order of the critical point is $k\geq2$ are
  not covered. This is because with this assumption any
  reconstruction, regardless of the degree of the corresponding
  polynomials, will attain an order of at least $k+1\geq3$. Therefore,
  the accuracy will be optimal regardless of the values of the weights
  $\omega_i$, taking into consideration that they always are a convex
  combination, namely, $\omega_0,\omega_1\geq0$ and
  $\omega_0+\omega_1=1$.
\end{remark}

\subsubsection*{Summary of the algorithm}

Input: $\bar{S}=\{f_{-1},f_0,f_1,f_2\}$, with $f_i=f(x_i)$ or $f_i=\frac{1}{h}\int_{x_{i-1/2}}^{x_{i+1/2}}f(x) \, \textrm{d}x$, and $\varepsilon>0$.
\begin{enumerate}
\item Compute the corresponding interpolating polynomials evaluated at $x_{1/2}$, which, both in case of reconstructions from point values and from cell averages, are given by \eqref{p0xhalf}. 
  
\item Compute the corresponding Jiang-Shu smoothness indicators $I_0$, $I_1$ and~$I_2$ (including the one considering the rightmost node) by~\eqref{eq:I0I1} and~\eqref{eq:I2}. 
\item Compute the auxiliary weights~$\tilde{\omega}_0$ and~$\tilde{\omega}_1$ from \eqref{eq:tw}. 
\item Define $\tau$ by \eqref{eq:taudef}. 
\item Compute the corrector weight~$\omega$ from \eqref{eq:cw}. 
\item Compute the corrected weights $\omega_0$ and $\omega_1$ from \eqref{correctedweights}. 
\item Obtain the OWENO reconstruction at $x_{1/2}$:
  $$p_2(x_{1/2})=\omega_0p_0(x_{1/2})+\omega_1p_1(x_{1/2}).$$
\end{enumerate}
Output: $\mathcal{R}(f_{-1}, f_0, f_1, f_2, \eps):=p_2(x_{1/2})$.

\subsubsection*{Comparison of the algorithm with the
  YC-WENO3 approach}

 Let us stress the key differences between the YC-WENO3 method and our proposal. Although at first sight both schemes might look similar, as they include squared divided differences including more than two nodes, the new method is not just an extension of YC-WENO3 including an additional downwind node. This additional node is used in the steps involving formulas
	\eqref{eq:I2} to
	\eqref{correctedweights}, in a manner so that the problem of order loss 
	at the critical points in the YC-WENO3 scheme is avoided.  
	This accuracy loss is due to the fact that the quotient $\sigma/I_i$ that appears in \eqref{eq:alphai} does not converge to zero as $h\to0^+$ when the point to which the stencil converges,
	$z$, satisfies $f'(z)=0$ and $f''(z)\neq0$ (namely, when $z$
	is a first-order critical point),
	since in that case there holds $\sigma=\bar\bigO(h^4)$ and
	$I_i=\bigO(h^4)$.
	 According to the results presented above, summarized in Theorem \ref{accuracy}, we have obtained a third-order WENO
	reconstruction procedure which, unlike the YC-WENO3 approach, attains
	the optimal third-order accuracy near critical points.
It is important to remark that this is in turn consistent with the
conclusion obtained in Section \ref{subsec:weno3}, in which it is
proven that there cannot exist a 3-point interpolator accounting for
discontinuities, while in turn maintaining the optimal third-order
accuracy near critical points, unless an artificially-scaled tuning
parameter is used.

\section{WENO schemes for systems of conservation laws}\label{sec:wenohyp}

In this section we discuss the incorporation of the novel third-order
WENO approach in the context of hyperbolic conservation laws. The
purpose is to prove that the resulting scheme depends on the same grid
points as a  standard  third-order WENO
reconstruction based on a three-point  stencil.

\subsection{Hyperbolic systems of conservation laws}

We will briefly describe in this section the equations and their discretization procedure.
We consider hyperbolic systems of~$\nu$ scalar 
conservation laws in $d$~space dimensions: 
\begin{align}\label{eq:hcl}
    \boldsymbol{u}_t+\sum_{i=1}^d \boldsymbol{f}^i(\boldsymbol{u})_{x_i}&= \boldsymbol{0},
     \quad (\boldsymbol{x},t)\in\Omega\times\mathbb{R}^+\subseteq\mathbb{R}^d\times\mathbb{R}^+, 
      \quad \boldsymbol{x} = (x_1,\ldots,x_d), 
\end{align} 
where $\boldsymbol{u}= \boldsymbol{u} ( \boldsymbol{x},t)\in\mathbb{R}^{\nu}$ is the sought solution, 
 $\boldsymbol{f}^i: \mathbb{R}^{\nu} \rightarrow\mathbb{R}^{\nu}$ are given flux density vectors, and 
\begin{align*}  
 \boldsymbol{u}=\begin{pmatrix} 
    u_1 \\
    \vdots \\
    u_{\nu}
    \end{pmatrix},\quad \boldsymbol{f}^i= \begin{pmatrix}       f^i_1 \\
      \vdots \\
      f^i_{\nu}
      \end{pmatrix}, \quad i=1, \dots, d;  \quad 
     \boldsymbol{f}=\begin{bmatrix} \boldsymbol{f}^1&\dots& \boldsymbol{f}^d
      \end{bmatrix}.
     \end{align*} 
System \eqref{eq:hcl} is complemented with  the initial condition 
\begin{align*}
u(\boldsymbol{x},0)=\boldsymbol{u}_0(\boldsymbol{x}),  \quad  \boldsymbol{x}\in\Omega, 
\end{align*} 
and prescribed boundary conditions. 

To describe the spatial discretization, we introduce a  Cartesian grid $\mathcal{G}$ formed by points (cell centers) $\smash{\boldsymbol{x}=\boldsymbol{x}_{j_1,\dots,j_d}=((j_1-\frac{1}{2})h,\dots,(j_d-\frac{1}{2})h)\in\mathcal{G}}$ for $h>0$. 
 In what follows, we use the index vector $\boldsymbol{j} = (j_1,
 \dots, j_d)$,  let $\boldsymbol{e}_i$ denote the $i$-th
 $d$-dimensional unit vector, and assume that $J$ is the set of all
 indices $\boldsymbol{j}$ for which point values of the solution are
 to be computed. We then advance  a semi-discrete scheme in which
 spatial derivatives are discretized first. The result is a system of ordinary
 differential equations whose numerical solution is iteratively
 updated in time. To do so, we first define 
 \begin{align*} 
  \boldsymbol{U} (t) := \bigl( \boldsymbol{u} ( \boldsymbol{x}_{\boldsymbol{j}}, t ) \bigr)_{\boldsymbol{j}\in J}. 
  \end{align*}   
To solve \eqref{eq:hcl}
we utilize  the Shu-Osher finite difference scheme \cite{shuosher88,shuosher89} with upwind spatial reconstructions of the flux function that are incorporated  into numerical flux vectors  $\smash{\boldsymbol{\hat f}}^{i}$ through a Donat-Marquina flux-splitting \cite{DonatMarquina96}. Thus, the contribution to the flux divergence in the coordinate~$x_i$ at point
 $\smash{\boldsymbol{x}=\boldsymbol{x}_{\boldsymbol{j}}}$  is given by 
\begin{align*} 
\boldsymbol{f}^i( \boldsymbol{U})_{x_i}(\boldsymbol{x}_{\boldsymbol{j}}, t) &  \approx\frac{1}{h} \Bigl( 
  \smash{\boldsymbol{\hat f}}^i_{\boldsymbol{j} + \frac{1}{2} \boldsymbol{e}_i} \bigl(\boldsymbol{U} ( t)\bigr)-
   \smash{\boldsymbol{\hat f}}^i_{ \boldsymbol{j} - \frac{1}{2} \boldsymbol{e}_i} \bigl( \boldsymbol{U}( t) \bigr) \Bigr). 
 \end{align*}
Then, WENO reconstructions \cite{JiangShu96} of order $2r+1$ are
considered, with special emphasis on the case we are interested in,
namely, $r=1$ (order $3$). 
To specify the time discretization, we write the semi-discrete scheme compactly as 
\begin{align*} 
\frac{\mathrm{d}} {\mathrm{d} t} \boldsymbol{U} (t) = \boldsymbol{\mathcal{L}} ( \boldsymbol{U} (t) ), 
\quad  \boldsymbol{\mathcal{L}} \bigl( \boldsymbol{U} (t) \bigr) =  \bigl( \mathcal{L}_{\boldsymbol{j}}  ( \boldsymbol{U} (t) ) \bigr)_{\boldsymbol{j} 
 \in J}, 
\end{align*} 
where we define 
\begin{align*} 
\mathcal{L}_{\boldsymbol{j}}  ( \boldsymbol{U} (t) ):= \frac{1}{h} \sum_{i=1}^d \Bigl( 
  \smash{\boldsymbol{\hat f}}^i_{\boldsymbol{j} + \frac{1}{2} \boldsymbol{e}_i} \bigl(\boldsymbol{U} ( t)\bigr)-
   \smash{\boldsymbol{\hat f}}^i_{ \boldsymbol{j} - \frac{1}{2} \boldsymbol{e}_i} \bigl( \boldsymbol{U}( t) \bigr) \Bigr)
\end{align*} 
(with suitable modifications for boundary points). 

For the time discretization, we use either the third-order
TVD Runge-Kutta  scheme proposed in \cite{shuosher89} or the
approximate Lax-Wendroff (henceforth, LWA) approach proposed in
\cite{ZorioBaezaMulet17}, which in turn is based on the original
Lax-Wendroff (henceforth, LW) approach proposed by Qiu and Shu in
\cite{QiuShu2003}. The choice for the time discretization will be
specified in each numerical experiment.
 
\subsection{Third-order WENO scheme}\label{subsec:dependence}

Although it may seem that the overall scheme for finite dimensional
conservation laws uses more points than the corresponding scheme for
classical WENO3 reconstructions, it is not the case, as we now show.

The semidiscrete scheme for a scalar one-dimensional law is
\begin{align}\label{eq:pep1}
  u_{i}'(t)&=-\frac{1}{h}\bigl(
    \hat f_{i+1/2}-\hat f_{i-1/2}
  \bigr),\\
  \label{hatfiph} 
  \hat f_{i+1/2}&=\hat f(u_{i-1}, u_{i}, u_{i+1},u_{i+2}),
\end{align}
so that the right-hand side of \eqref{eq:pep1} depends on
approximations $u_{j}(t)\approx u(x_{j}, t)$ at a 5-point
stencil $j=i-2,\dots,i+2$. An ODE solver, such as the third-order
TVD Runge-Kutta scheme 
proposed in \cite{shuosher89}, is applied to \eqref{eq:pep1} to obtain
the final time-space accurate scheme.

If the reconstruction \eqref{hatfiph} associated with the 
   cell interface~$x_{i+1/2}$ is sought and we define the interval $I(a,b) :=[\min\{a, b\},\max\{a, b\}]$,  then we determine 
    for $j \in \{ i-1,i, i+1, i+1\}$ the quantities 
   \begin{align*}
    f_j^{i+1/2,+} &:=
    \begin{cases}
      f(u_j) & \text{if $f'(u) > 0$ for all $u\in I(u_{i},u_{i+1})$,} \\
      0 & \text{if $f'(u) < 0$  for all $u\in I(u_{i},u_{i+1})$,} \\
      f(u_j)+\alpha_{i+1/2}u_j&\text{otherwise,}
    \end{cases}  \\
    f_j^{i+1/2,-} &:=
    \begin{cases}
      0 & \text{if $f'(u) > 0$ for all $u\in I(u_{i},u_{i+1})$,} \\
      f(u_j) & \text{if $f'(u) < 0$ for all $u\in I(u_{i},u_{i+1})$,} \\
      f(u_j)-\alpha_{i+1/2}u_j& \text{otherwise,} 
    \end{cases}
  \end{align*}
     where 
  \begin{align*} 
  \alpha_{i+1/2}:=\max_{u\in I(u_{i}, u_{i+1})} \bigl|f'(u) \bigr|.
   \end{align*} 

The precise formulation (see \cite{shuosher89}) for attaining third-order
accuracy (the 
maximum for semidiscrete stability being three for this 5-point
stencil, cf.\ \cite{maxorder}) for the usual WENO3
reconstructions consists in using a local 
flux splitting  $f(u)= f^+(u) + f^-(u)$, 
such that
$\smash{\pm(f^{\pm}(u))' \geq 0}$, in the interval $I(u_{i},u_{i+1}]$ determined by $u_{i}$ 
  and $u_{i+1}$, which is defined as 
  \begin{align}  \label{reconstr-standard}  \begin{split} 
   \hat{f}_{i+1/2}  & := \mathcal{R}^+ 
    \bigl( f^{i+1/2,+}_{i-1}, f^{i+1/2,+}_i, f^{i,+1/2,+}_{i+1} \bigr) \\ & \quad  + \mathcal{R}^- 
        \bigl( f^{i+1/2,-}_{i}, f^{i+1/2,-}_{i+1}, f^{i+1/2,-}_{i+2} \bigr),   \end{split} 
        \end{align} 
  where $\mathcal{R}^{+}$ is a right-biased cell-averages
  reconstruction and is the right-biased cell-averages reconstruction
  given by  $\mathcal{R}^{-}(a, b, c)=\mathcal{R}^{+}(c, b, a)$. In contrast, the flux splitting and reconstruction used herein are defined as follows. Instead of using \eqref{reconstr-standard}, we propose to 
    define the flux value    $\smash{\hat f_{i+1/2}}$ 
  by our optimal-order reconstruction    
        $\smash{\mathcal{R}^{\pm,\text{opt}}}$  that depends on
  the  four-point  stencil, such that  
  \begin{align*}
    \hat f_{i+1/2} &=\hat f(u_{i-1}, u_{i}, u_{i+1},u_{i+2})\\
    &=
    \mathcal{R}^{+,\text{opt}} \bigl(f^{i+1/2,+}_{i-1}, f^{i+1/2,+}_{i},f^{i+1/2,+}_{i+1},f^{i+1/2,+}_{i+2} \bigr) \\
    & \quad +
    \mathcal{R}^{-,\text{opt}} \bigl(f^{i+1/2,-}_{i-1},f^{i+1/2,-}_{i},f^{i+1/2,-}_{i+1}  ,f^{i+1/2,-}_{i+2} \bigr).
  \end{align*}

Systems of conservation laws are dealt by the application of the former scheme
to local characteristic fields, obtained by a double linearization
\cite{DonatMarquina96}. The extension to multidimensional Cartesian
grids follows by working dimension by dimension.

It is also important to remark that, as pointed out in \cite{maxorder},
 the maximum order of accuracy that a stable semidiscrete scheme with a numerical domain of dependence of five points can attain is three. This explains why we
  are not using the additional fourth node (used only for the
  smoothness analysis) to increase by a unit the
  order of the reconstructions. Otherwise, and in this particular
  case, the reconstructions would be centered, thus our scheme
  would lose the upwind features, and ultimately the stability
  properties.
  
  \section{Numerical experiments}\label{sec:num_exp}

\subsection{Accuracy tests with algebraic problems}

We perform some tests in order to verify the accuracy properties of
the scheme in presence of critical points. To this end, we use the multiple-precision library MPFR \cite{MPFR} through its C++ wrapper \cite{Holoborodko}, using a precision of $3322$ bits ($\approx1000$ digits) and taking in all cases $\varepsilon=10^{-10^6}$.

Let us consider the family of functions $f_k:\mathbb{R}\to\mathbb{R}$, $k \in \{ 0, 1\}$, given by
$f_k(x)=x^{k+1}\mathrm{e}^x$. Then $f_k$~has a smooth extremum at~$x=0$ of order~$k$. 
In this case the error is given by $E_{k,n}=|P_n(0)-f_k(0)|$, where $P_n$ denotes the corresponding reconstruction at $x_{1/2}=0$, with the grid $x_i=(i- 1/2)h$, $i\in\{-1,\ldots,1+s\}$ ($s=0$ for the traditional JS-WENO and YC-WENO schemes and $s=1$ for the proposed optimal WENO schemes, in which an additional node is considered), with $h=1/n$ for $n\in\N$, when pointwise values are taken, namely, $f_{k,i}=f_k(x_i)$ and reconstructions from pointwise values to pointwise values are performed. We also 
 present in the tables the same setup when cell average values are taken instead: 
 \begin{align*} 
 f_{k,i}=\int_{x_i- h/2}^{x_i+ h / 2}f(x)\, \mathrm{d}x, 
 \end{align*} 
 by reconstructing pointwise values from cell average values. In all cases, the tables show the corresponding average
 orders,\\ $O_k=\frac{1}{80}\sum_{j=1}^{80}o_{k,j}$, where
 $o_{k,j}=\log_2(E_{k,n_{j-1}}/E_{k,n_j}))$, with $n_j=5\cdot2^j$,
 $j\in\{0,\ldots,80\}$. 
 
We consider alternatively  the family of functions $g_k:\mathbb{R}\to\mathbb{R}$, $k \in \{ 0, 1\}$ given by
\[g_k(x)=
\begin{cases}
  x^{2k}\mathrm{e}^x & \text{for $x\leq0$,} \\
  \mathrm{e}^{x+1} & \text{for $x>0$.} 
\end{cases}
\] 
Then $g_k$ has a discontinuity at $x=0$ with a left smooth extremum of order~$k$ for $k \in \{ 0,1\}$. We test the accuracy of the methods with the same parameters as above, where, in order to  emphasize the behaviour of our optimal scheme at discontinuities, in this case we change the location of the discontinuity by considering a grid of the form $x_i=(i-\frac{1}{2}+\theta)h$, $i\in\{-1,\ldots,1+s\}$, for $\theta\in\{0,1\}$. Since $x_{1/2}=\theta h$, the error is thus now given by $|P(\theta h)-g(\theta h)|$.

The results involving the different combinations of the proposed values for $k$ in the case of $f_k$ and for $\theta$ and $k$ in the case of $g_k$ are shown in Table \ref{alg} for the traditional JS-WENO3 and YC-WENO3 schemes as well as the optimal WENO approach presented herein.

\begin{table}
  \begin{center} 
  \begin{tabular}{|c|c|c|c|c|c|c|c|c|}
    \hline
    & \multirow{2}{*}{$\theta$} & \multirow{2}{*}{$k$} & \multicolumn{2}{|c|}{JSWENO3} & \multicolumn{2}{|c|}{YCWENO3} & \multicolumn{2}{|c|}{OWENO3} \\
    \cline{4-9}
    & & & Point & Cell & Point & Cell & Point & Cell \\
    \hline
    \multirow{2}{*}{$f_k$} & --- & $0$ & 3.00 & 3.00 & 2.98 & 2.98 & 3.00 & 3.00 \\
    & --- & $1$ & 2.00 & 2.00 & 2.00 & 2.00 & \textbf{3.01} & \textbf{3.01} \\
    \hline
    \multirow{4}{*}{$g_k$} & $0$ & $0$ & 1.97 & 1.93 & 1.98 & 1.98 & 2.00 & 2.00 \\
    & $0$ & $1$ & 1.99 & 1.99 & 1.99 & 1.99 & 1.96 & 1.95 \\
    & $1$ & $0$ & 2.00 & 2.00 & 2.00 & 2.00 & 1.99 & 2.00 \\
    & $1$ & $1$ & 2.00 & 2.00 & 2.00 & 2.00 & 2.00 & 2.00 \\
    \hline
  \end{tabular} \end{center}  \smallskip 
  \caption{Numerical order for third-order schemes, functions with
    smooth extrema.}
  \label{alg}
\end{table}

We discuss row by row the results obtained in Table~\ref{alg}. The first two rows containing data stand for the function $f_k$, which is a smooth function with a critical point of order~$k$. Therefore, the optimal order is $3$. We can see that when the critical point has order zero, namely, $k=0$, all the schemes attain the optimal accuracy. However, differences arise when $k=1$. In this case, the first-order critical point affects the traditional WENO schemes decreasing its accuracy in one unit, whereas the optimal WENO approach keeps the optimal third-order accuracy.

As for the function $g_k$, we can conclude that regardless of the
position of the discontinuity with respect to the stencil $S$ and the
order of the critical point, all the
schemes, both the traditional ones and the optimal ones, attain the
suboptimal second-order accuracy, avoiding the error of magnitude
$\bar{\bigO}(1)$ associated to the substencil containing the
discontinuity. This is the best order of accuracy that can be obtained
near  a discontinuity by shock-capturing methods based on three- or four-point stencils.

\subsection{Conservation law experiments}

In this section we present some experiments involving numerical solutions of hyperbolic conservation laws. We discretize them in time by the approximate Lax-Wendroff approach matching the spatial order proposed in \cite{ZorioBaezaMulet17}, unless we indicate the contrary, in whose case the third-order
TVD Runge-Kutta scheme \cite{shuosher89} will be used. Also, since in
this case we work with double precision, the $\varepsilon$ parameter
is chosen as $\varepsilon=10^{-100}$. The flux splitting used is
Donat-Marquina \cite{DonatMarquina96} for the problems with weak
solutions and Local Lax-Friedrichs for the
problems with smooth solutions (unless all the characteristics move to
the same direction, in whose case we simply use the corresponding
left/right-biased upwind reconstructions). In all cases, and also
unless we state the contrary, the CFL used for the 1D experiments is
$0.5$ and $0.4$ for the 2D experiments. The
reason for the choice of these CFL values is for uniformity reasons,
combined with the fact that in
some problems with complex structures or interactions between
discontinuities, the fifth-order WENO method combined with the flux
splitting used \cite{DonatMarquina96} can develop some
oscillations (such as in the double Mach reflection problem, in Example
5) or even fail (such as in the blast wave problem, in Example 4) if
larger CFL values are used. It must be
pointed out that these issues have not been observed on any of the
third-order methods for larger CFL values.

\subsubsection*{Example 1: Linear advection equation}

We consider the linear advection equation with the following domain, boundary condition and initial condition:
\begin{gather*} 
u_t+f(u)_x=0,\quad\Omega=(-1,1),\quad u(-1, t)=u(1, t), \\
f(u)=u,\quad u_0(x)=0.25+0.5\sin(\pi x), 
\end{gather*} 
whose exact solution is $u(x,t)=0.25+0.5\sin(\pi(x-t))$, with critical points located at $x=t+m+1/2$, $m\in\Z$. We run several simulations with final time $T=1$, resolutions of $n$ points, that is, with a grid spacing of $h=2/n$, using the classical JS-WENO schemes, YC-WENO and our OWENO3 scheme, both with the $\|\cdot\|_1$ and $\|\cdot\|_{\infty}$ errors. Since the characteristics move to the right, we use left-biased reconstructions. The results are shown in Table~\ref{advlin_o3}. From the table   it can be appreciated that an accuracy loss is produced in the case of the traditional schemes, whereas the optimal third-order accuracy is solidly kept by the novel scheme.

\begin{table}
  \setlength\tabcolsep{2.5pt}
  \centering  \footnotesize
  \begin{tabular}{|c|c|c|c|c|c|c|c|c|c|c|c|c|}   
    \hline
    & \multicolumn{4}{c|}{JSWENO3} & \multicolumn{4}{c|}{YCWENO3} & \multicolumn{4}{c|}{OWENO3} \\
    \hline
    & \multicolumn{2}{c|}{$\|\cdot\|_1$} & \multicolumn{2}{c|}{$\|\cdot\|_{\infty}$} & \multicolumn{2}{c|}{$\|\cdot\|_1$} & \multicolumn{2}{c|}{$\|\cdot\|_{\infty}$} & \multicolumn{2}{c|}{$\|\cdot\|_1$} & \multicolumn{2}{c|}{$\|\cdot\|_{\infty}$} \\
    \hline
    $n$ & Err. & $\bigO$ & Err. & $\bigO$ & Err. & $\bigO$ & Err. & $\bigO$ & Err. & $\bigO$ & Err. & $\bigO$ \\
    \hline
    40 & 8.52e-03 & --- & 2.56e-02 & --- & 6.67e-03 & --- & 2.11e-02 & --- & 1.87e-04 & --- & 3.08e-04 & --- \\
    80 & 2.10e-03 & 2.02 & 1.00e-02 & 1.36 & 1.46e-03 & 2.19 & 7.90e-03 & 1.42 & 2.31e-05 & 3.02 & 3.66e-05 & 3.07 \\
    160 & 4.86e-04 & 2.11 & 3.81e-03 & 1.39 & 3.19e-04 & 2.20 & 2.87e-03 & 1.46 & 2.86e-06 & 3.01 & 4.50e-06 & 3.02 \\
    320 & 1.10e-04 & 2.15 & 1.43e-03 & 1.42 & 6.45e-05 & 2.31 & 1.02e-03 & 1.50 & 3.56e-07 & 3.01 & 5.60e-07 & 3.01 \\
    640 & 2.45e-05 & 2.16 & 5.28e-04 & 1.43 & 1.32e-05 & 2.29 & 3.54e-04 & 1.52 & 4.44e-08 & 3.00 & 6.98e-08 & 3.00 \\
    1280 & 5.42e-06 & 2.18 & 1.94e-04 & 1.45 & 2.61e-06 & 2.34 & 1.21e-04 & 1.55 & 5.55e-09 & 3.00 & 8.72e-09 & 3.00 \\
    2560 & 1.19e-06 & 2.19 & 7.06e-05 & 1.46 & 5.05e-07 & 2.37 & 4.10e-05 & 1.56 & 6.93e-10 & 3.00 & 1.09e-09 & 3.00 \\
    5120 & 2.57e-07 & 2.21 & 2.56e-05 & 1.46 & 9.69e-08 & 2.38 & 1.37e-05 & 1.58 & 8.67e-11 & 3.00 & 1.36e-10 & 3.00 \\
    10240 & 5.54e-08 & 2.21 & 9.22e-06 & 1.47 & 1.84e-08 & 2.40 & 4.53e-06 & 1.60 & 1.08e-11 & 3.00 & 1.71e-11 & 2.99 \\
    20480 & 1.19e-08 & 2.22 & 3.31e-06 & 1.48 & 3.44e-09 & 2.42 & 1.49e-06 & 1.61 & 1.43e-12 & 2.92 & 2.38e-12 & 2.84 \\
    \hline
  \end{tabular}
  \smallskip 
  
  \caption{Example 1: linear advection equation, third-order schemes.}
  \label{advlin_o3}
\end{table}

\begin{table}
  \setlength\tabcolsep{2.5pt}
  \centering \footnotesize
  \begin{tabular}{|c|c|c|c|c|c|c|c|c|c|c|c|c|}
    \hline
    & \multicolumn{4}{c|}{JSWENO3} & \multicolumn{4}{c|}{YCWENO3} & \multicolumn{4}{c|}{OWENO3} \\
    \hline
    & \multicolumn{2}{c|}{$\|\cdot\|_1$} & \multicolumn{2}{c|}{$\|\cdot\|_{\infty}$} & \multicolumn{2}{c|}{$\|\cdot\|_1$} & \multicolumn{2}{c|}{$\|\cdot\|_{\infty}$} & \multicolumn{2}{c|}{$\|\cdot\|_1$} & \multicolumn{2}{c|}{$\|\cdot\|_{\infty}$} \\
    \hline
    $n$ & Err. & $\bigO$ & Err. & $\bigO$ & Err. & $\bigO$ & Err. & $\bigO$ & Err. & $\bigO$ & Err. & $\bigO$ \\
    \hline
    40 & 1.77e-03 & --- & 1.11e-02 & --- & 1.62e-03 & --- & 9.85e-03 & --- & 1.70e-04 & --- & 1.11e-03 & --- \\
    80 & 4.77e-04 & 1.89 & 4.17e-03 & 1.41 & 4.21e-04 & 1.95 & 3.57e-03 & 1.46 & 2.24e-05 & 2.92 & 1.81e-04 & 2.62 \\
    160 & 1.18e-04 & 2.02 & 1.62e-03 & 1.36 & 9.80e-05 & 2.10 & 1.32e-03 & 1.43 & 2.75e-06 & 3.03 & 2.27e-05 & 3.00 \\
    320 & 2.91e-05 & 2.02 & 6.21e-04 & 1.38 & 2.24e-05 & 2.13 & 4.94e-04 & 1.42 & 3.37e-07 & 3.03 & 2.77e-06 & 3.03 \\
    640 & 7.01e-06 & 2.06 & 2.36e-04 & 1.40 & 5.05e-06 & 2.15 & 1.79e-04 & 1.46 & 4.16e-08 & 3.02 & 3.41e-07 & 3.02 \\
    1280 & 1.64e-06 & 2.10 & 8.84e-05 & 1.42 & 1.11e-06 & 2.18 & 6.36e-05 & 1.50 & 5.17e-09 & 3.01 & 4.23e-08 & 3.01 \\
    2560 & 3.83e-07 & 2.10 & 3.28e-05 & 1.43 & 2.47e-07 & 2.17 & 2.21e-05 & 1.52 & 6.44e-10 & 3.00 & 5.26e-09 & 3.01 \\
    5120 & 8.85e-08 & 2.11 & 1.20e-05 & 1.45 & 5.45e-08 & 2.18 & 8.84e-06 & 1.32 & 8.04e-11 & 3.00 & 6.56e-10 & 3.00 \\
    10240 & 2.04e-08 & 2.11 & 5.42e-06 & 1.15 & 1.21e-08 & 2.17 & 4.41e-06 & 1.00 & 1.00e-11 & 3.00 & 8.19e-11 & 3.00 \\
    20480 & 4.70e-09 & 2.12 & 2.73e-06 & 0.99 & 2.71e-09 & 2.16 & 2.21e-06 & 1.00 & 1.25e-12 & 3.00 & 1.02e-11 & 3.00 \\
    \hline
  \end{tabular} 
  \smallskip 
  \caption{Example 2a: Burgers equation, third-order schemes.}
  \label{burgers_o3}
\end{table}

\subsubsection*{Examples 2a and 2b: Burgers equation}

\begin{figure}[t] 
  \centering
  \includegraphics{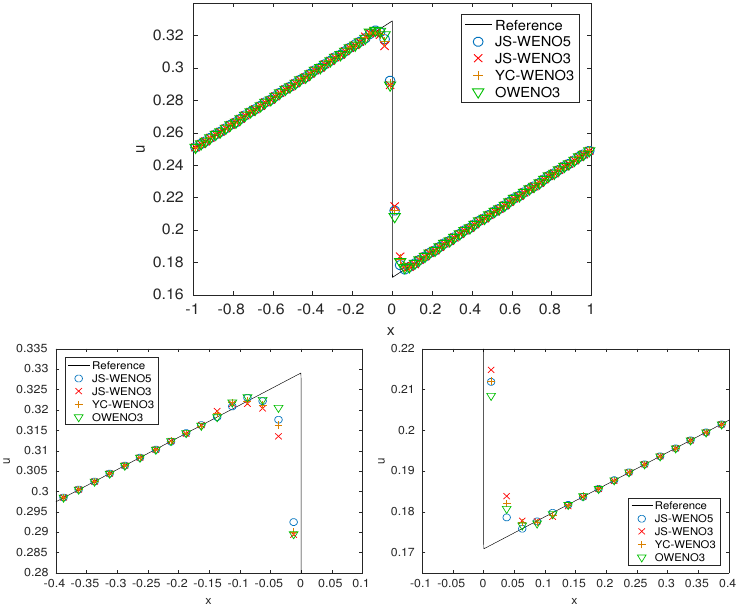} 
    \caption{Example 2b (Burgers equation, discontinuous solution at $T=12$): 
   third-order schemes.}
  \label{burgersdisc_o3}
\end{figure}

We now consider Burgers equation with the following setup involving the domain, boundary conditions and initial condition:
\begin{gather} \label{burgersprob} \begin{split} 
  u_t+f(u)_x=0,\quad\Omega=(-1,1),\quad u(-1, t)=u(1, t), \\
 f(u)=u^2/2, \quad 
u_0(x)=0.25+0.5\sin(\pi x).
 \end{split} 
\end{gather}
In this case, $f(u_0(x))$ has first-order smooth extrema at $x=-1/2$ and at $x=1/2$. 
In  Example~2a, we consider  the solution of \eqref{burgersprob} at  $T=0.3$, 
when it is still smooth, whose results are shown in Table \ref{burgers_o3}, while in Example~2b we set $T=12$, when the solution of \eqref{burgersprob} has become discontinuous, shown in Figure \ref{burgersdisc_o3}, in which are also compared against the results obtained by the widely used JS-WENO5 schemes.

From Table \ref{burgers_o3} one can see that again, as in the linear advection case, the presence of first-order critical points makes the accuracy of the traditional schemes decay to orders lower than three, while the third-order accuracy is still kept by the optimal third-order scheme. As for the discontinuous case, we can see in Figure \ref{burgersdisc_o3} that the optimal third-order scheme has a much higher resolution than the traditional 
third-order schemes, especially near the discontinuity and, moreover, it is similar to the resolution presented by the fifth-order scheme.

\subsubsection*{Examples 3a and 3b: Shu-Osher problem}

The 1D Euler equations for gas dynamics are given by 
 $\boldsymbol{u} = ( \rho, \rho v, E)^{\mathrm{T}}$ and $\boldsymbol{f} (\boldsymbol{u}) = 
  \boldsymbol{f}^1  (\boldsymbol{u}) = (\rho v, 
        p+\rho v^2, 
        v(E+p))^{\mathrm{T}}$,  where $\rho$ is the density, $v$ is the velocity, $E$ is the specific energy of the system and $p$ is the pressure, given by the equation of state
$p=(\gamma-1)(E-\rho v^2/2)$, 
where $\gamma$ is the adiabatic constant that will be taken as $\gamma =1.4$. 
We now consider the interaction with a Mach 3 shock and a sine wave. The spatial domain is now given by $\Omega:=(-5,5)$  with the initial condition 
\begin{align*} 
(\rho,v,p) (x, 0) = 
\begin{cases}
  (27/7, 
    4\sqrt{35} / 9, 
    31 / 3 ) & \text{if  $x\leq-4$,}  \\ 
   (1+\sin(5x)/5, 0, 1 ) 
    & \text{if $x>-4$,} 
\end{cases} \end{align*} 
with left inflow and right outflow boundary conditions.

\begin{figure}[t] 
  \centering
  \includegraphics{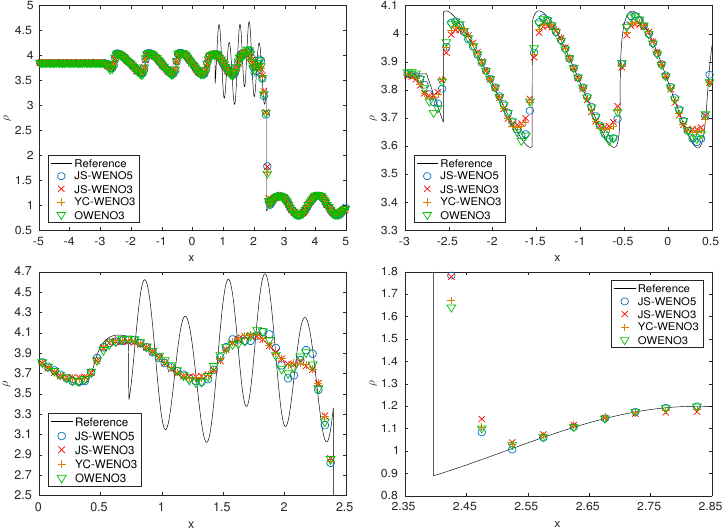} 
    \caption{Example 3a: Shu-Osher problem. $T=1.8$. $n=200$.}
  \label{shuosher200}
\end{figure}

\begin{figure}[t] 
  \centering
  \includegraphics{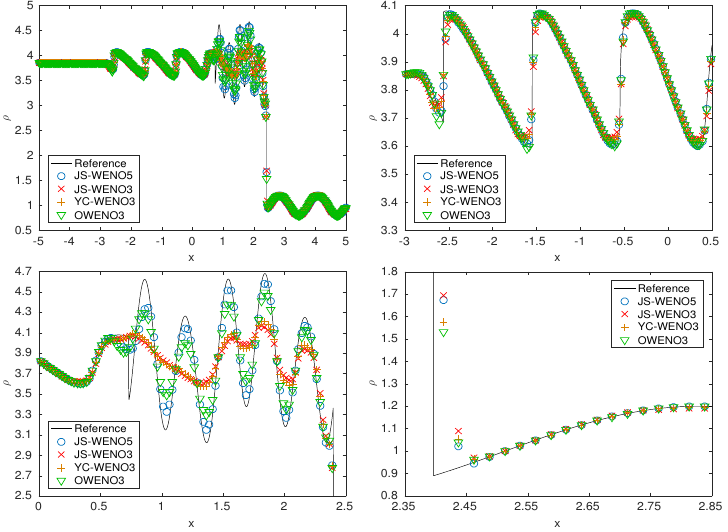} 
  \caption{Example 3a: Shu-Osher problem. $T=1.8$. $n=400$.}
  \label{shuosher400}
\end{figure}

We run the simulation until $T=1.8$ and compare the schemes against a reference solution computed with a resolution of $n=16000$ cells. Figures \ref{shuosher200} and \ref{shuosher400} are associated to the third-order schemes and JS-WENO5 with resolutions of $n=200$ and $n=400$ points, respectively, showing the corresponding density fields.

\begin{figure}[t] 
  \centering
  \includegraphics[width=0.96\textwidth]{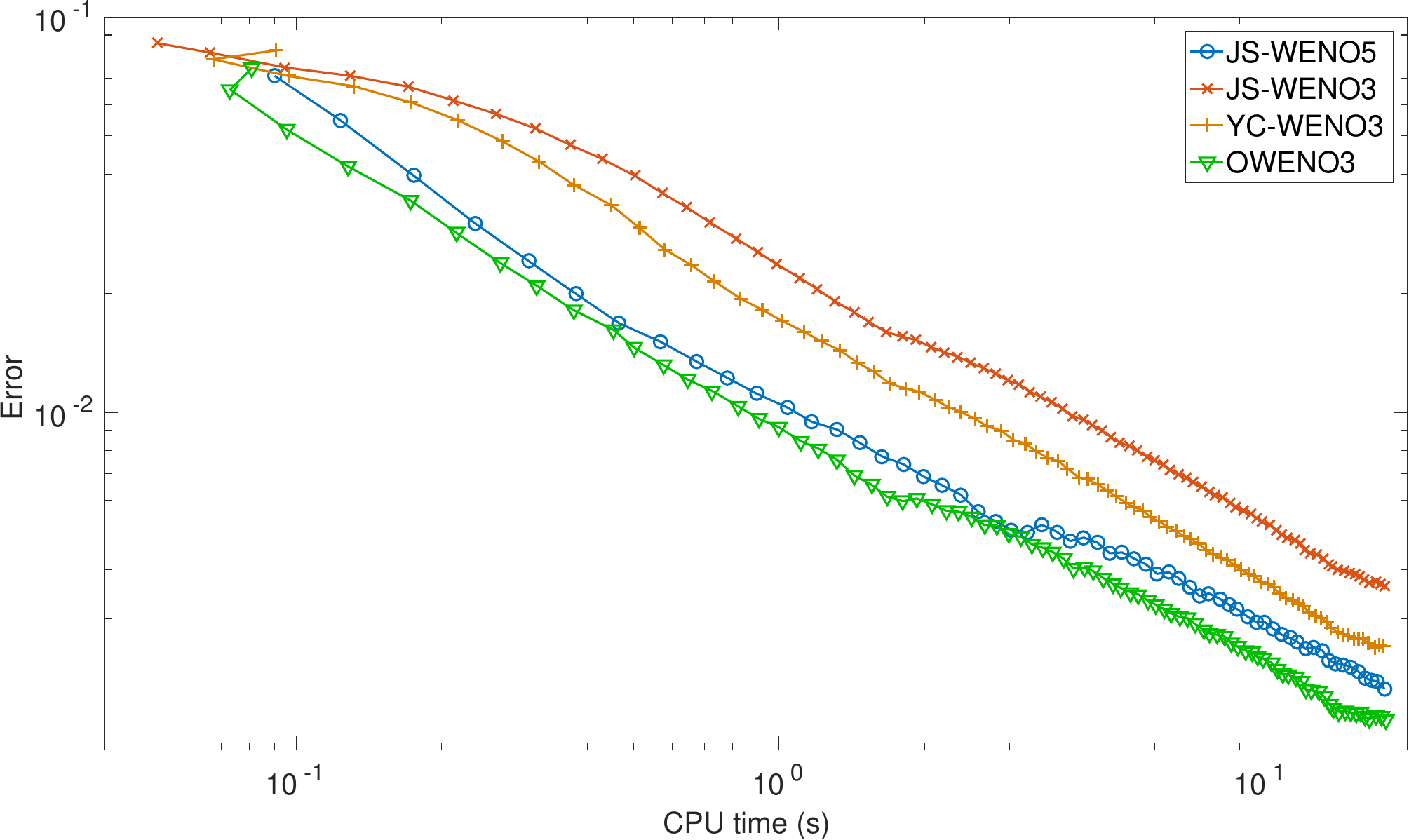} 
  \caption{Example 3b: Ratio Error/CPU comparison for Shu-Osher problem.}
  \label{shuosher_cpu}
\end{figure}

The third-order optimal schemes show again a much better resolution than their traditional counterparts, especially observed in the resolution of $n=400$ cells. Moreover, they have a similar resolution than the JS-WENO5 scheme, and at lower computational cost. In order to support the latter statement, we next present an efficiency comparison involving the ratio error $\|\cdot\|_1$ / CPU time, which can be seen in Figure~\ref{shuosher_cpu}.
The proposed third-order scheme shows a better performance than its traditional counterparts. Moreover, it is also more efficient than the JS-WENO5 scheme in this case.

\subsubsection*{Examples 4a and 4b: Blast wave problem}

Continuing with the 1D Euler equations, let us now simulate the interaction of two blast waves \cite{Colella} by using the
following initial data
$$u(x,0)=\begin{cases}
    u_{\mathrm{L}} & \text{if }0<x<0.1,\\
    u_{\mathrm{M}} & \text{if }0.1<x<0.9,\\
    u_{\mathrm{R}} & \text{if }0.9<x<1,
    \end{cases}$$
where $\rho_{\mathrm{L}}=\rho_{\mathrm{M}}=\rho_{\mathrm{R}}=1$, $v_{\mathrm{L}}=v_{\mathrm{M}}=v_{\mathrm{R}}=0$,
$p_{\mathrm{L}}=10^3,p_{\mathrm{M}}=10^{-2},p_{\mathrm{R}}=10^2$. We set reflecting boundary conditions
at $x=0$ and $x=1$, simulating a solid wall at both sides. This
problem involves multiple reflections of shocks and rarefactions off
the walls and many interactions of waves inside the domain.

\begin{figure}[t] 
  \centering
   \includegraphics{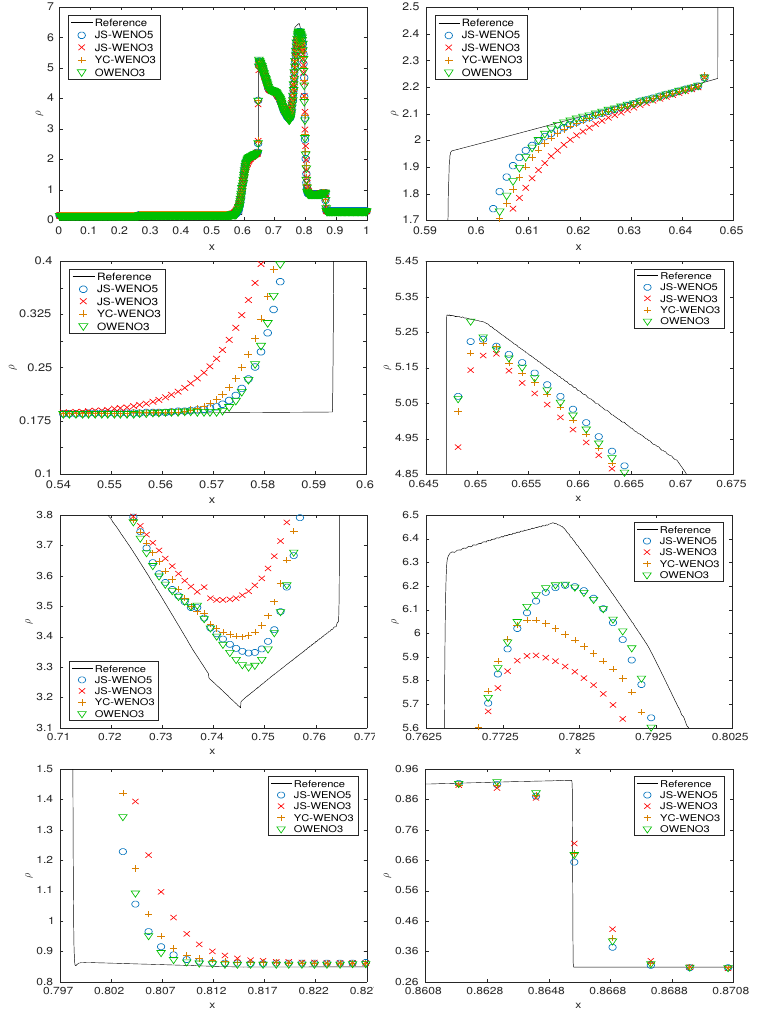}
   \caption{Example 4a: Blast wave problem. $T=0.038$.}
  \label{blast}
\end{figure}

The results are shown in Figure \ref{blast} for the density field at a resolution of $n=800$ cells, in which all the third-order schemes involved in this paper are used, being in turn compared with the JS-WENO5 scheme. The resolution used for the reference solution is $n=100000$ cells.

\begin{figure}[t] 
  \centering
  \includegraphics[width=0.96\textwidth]{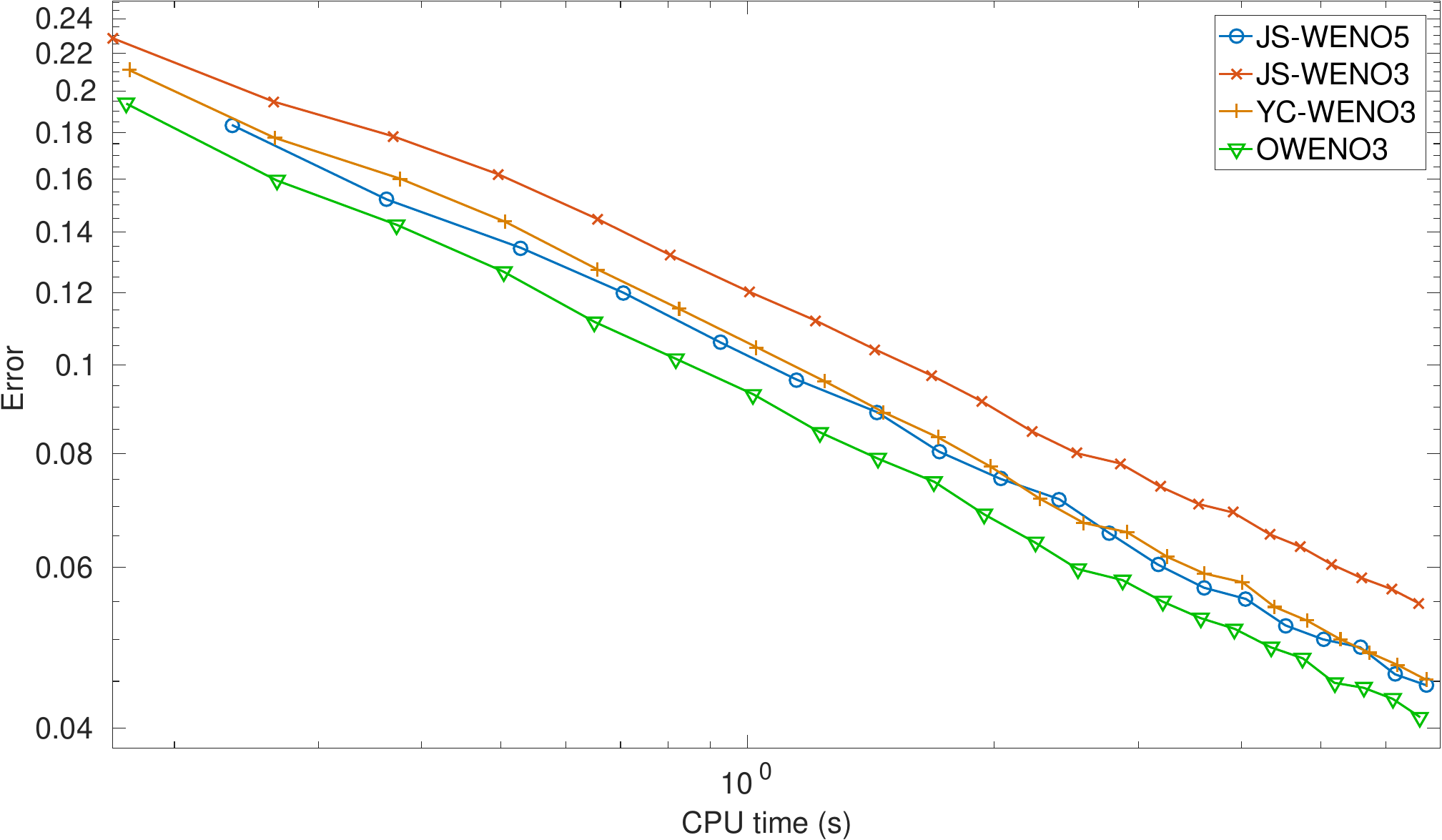}  
  \caption{Example 4b: Ratio Error/CPU comparison for blast wave problem.}
  \label{blast_cpu}
\end{figure}

As the results show, the third-order optimal scheme has at some regions a higher resolution than even the fifth-order scheme. Finally, Figure \ref{blast_cpu} shows an efficiency comparison between all the involved schemes, where, for the sake of performing a fair comparison, all the schemes have been equipped with the third-order
TVD Runge-Kutta scheme \cite{shuosher89}. In this case, the optimal third-order scheme is still more efficient than the fifth-order scheme.

\subsubsection*{Examples 5a and 5b: Double Mach reflection}

The equations that will be considered in this section are the two-dimensional Euler equations for inviscid gas dynamics
given by
\begin{align}\label{hcl}
    \boldsymbol{u}_t+\sum_{i=1}^d \boldsymbol{f}^i(\boldsymbol{u})_{x_i}&= \boldsymbol{0},
     \quad (\boldsymbol{x},t)\in\Omega\times\mathbb{R}^+\subseteq\mathbb{R}^d\times\mathbb{R}^+, 
      \quad \boldsymbol{x} = (x_1,\ldots,x_d), 
\end{align}
by taking in \eqref{hcl} $m=4$ and $d=2$, where setting $x=x_1$ and  $y=x_2$,  
 we have
\begin{align*} 
& \boldsymbol{u} =  \begin{pmatrix} 
    \rho \\
    \rho v^x \\
    \rho v^y \\
    E \end{pmatrix}, \quad 
  \boldsymbol{f}^1 (\boldsymbol{u})= \begin{pmatrix}  
    \rho v^x \\
    p+\rho (v^x)^2 \\
    \rho v^xv^y \\
    v^x(E+p) 
  \end{pmatrix},  \quad   \boldsymbol{f}^2 (\boldsymbol{u})= \begin{pmatrix} 
    \rho v^y \\
    \rho v^xv^y \\
    p+\rho (v^y)^2 \\
    v^y(E+p) 
  \end{pmatrix}. 
\end{align*}
Here  $\rho$ is the density, $(v^x, v^y)$  is the velocity, $E$ is the specific energy, and  $p$  is the pressure  that  is given by the equation of state $p=(\gamma-1)(E- \rho((v^x)^2+(v^y)^2)/2)$,  where
 the adiabatic constant is again chosen as  $\gamma =1.4$.

This experiment uses these equations to model a vertical right-going Mach 10 shock colliding with an equilateral triangle. By symmetry, this is equivalent to a collision with a ramp with a slope of $30^{\circ}$ with respect to the horizontal line.

\begin{figure}[t] 
  \centering
  \begin{tabular}{cc}
    \includegraphics[width=0.45\textwidth]{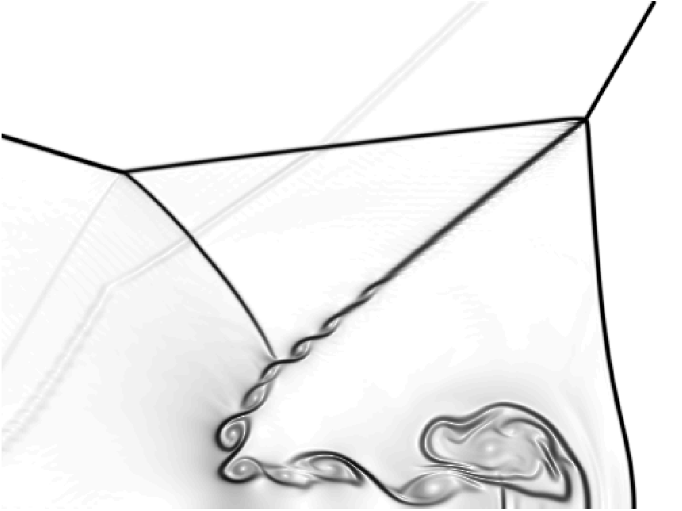} & \includegraphics[width=0.45\textwidth]{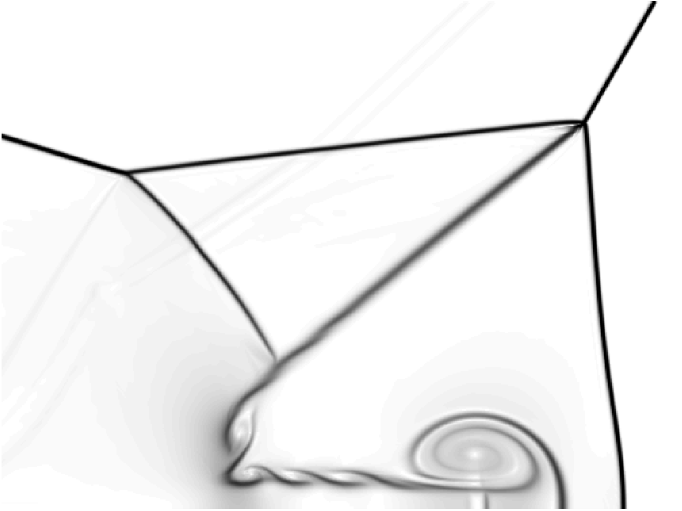} \\
    JS-WENO5 & JS-WENO3 \\
    \includegraphics[width=0.45\textwidth]{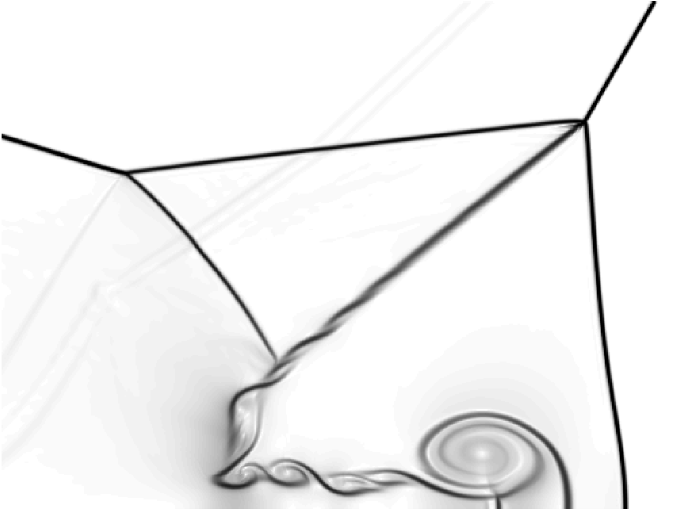} & \includegraphics[width=0.45\textwidth]{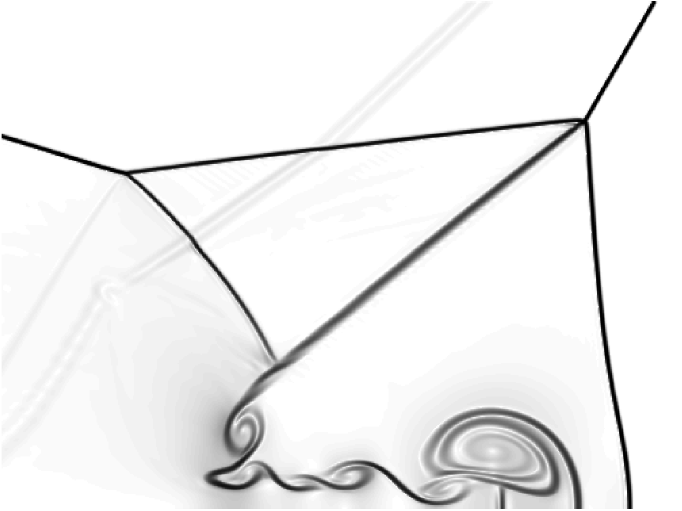} \\
    YC-WENO3 & OWENO3
  \end{tabular}
  \caption{Example 5a: Double Mach reflection,
    $2560\times640$. $T=0.2$. Schlieren plot of the density field.}
  \label{dmr640}
\end{figure}

For the sake of simplicity, in \cite{Colella} the equivalent problem is considered in a rectangle, consisting in a rotated shock, whose vertical angle is 
$30^{\circ}$. The domain is the rectangle $\Omega=[0,4]\times[0,1]$,
and the initial conditions are
\begin{gather*}  (\rho,v^x,v^y,E) (x,y,0)=\begin{cases}
  \boldsymbol{c}_1= (\rho_1,v_1^x,v_1^y,E_1)    & \text{if $y\leq\frac{1}{4}+\tan(\frac{\pi}{6})x$,} \\
  \boldsymbol{c}_2=  (\rho_2,v_2^x,v_2^y,E_2)     & \text{if $\frac{1}{4}+\tan(\frac{\pi}{6})x$,}  
\end{cases} \\
    \boldsymbol{c}_1 =
    \bigl(8,8.25\cos(\pi/6),-8.25\sin(\pi/6),563.5\bigr), \quad 
    \boldsymbol{c}_2=    (1.4,0,0,2.5).
  \end{gather*} 
We impose inflow boundary conditions, with value $\boldsymbol{c}_1$, at the left side, $\{0\}\times[0,1]$, outflow boundary conditions both at $[0,\frac{1}{4}]\times\{0\}$ and $\{4\}\times[0,1]$, reflecting boundary conditions at  $]\frac{1}{4},4]\times\{0\}$ and inflow boundary conditions at the upper side, $[0,4]\times\{1\}$, which mimics the shock at its actual traveling speed:
\begin{align*} 
 (\rho,v^x,v^y,E) (x,1,t)=\begin{cases}
  \boldsymbol{c}_1 & \text{if $x\leq\frac{1}{4}+\frac{1+20t}{\sqrt{3}}$,}  \\[2mm]
  \boldsymbol{c}_2 & \text{if $x>\frac{1}{4}+\frac{1+20t}{\sqrt{3}}$.}  
\end{cases} \end{align*} 
We run different simulations until $T=0.2$ at a resolution of $2560\times640$ points, shown in Figure \ref{dmr640}, with $\textnormal{CFL}=0.4$ and involving the classical JS-WENO5 scheme and the third-order schemes considered along this paper.


In this case, we can see that in both resolutions, both the YC-WENO3
scheme and the OWENO3 scheme have a higher resolution than
the JS-WENO3 scheme, in which the discontinuities and the non-smooth
features such as turbulence and vorticity are more smeared. On
the other hand, the resolution shown by the former schemes is still
remarkably lower than the JS-WENO5 scheme. This is probably due to the
nature of this problem, which has no solution for the inviscid 2D
Euler equations, since more and more turbulent structures
appear at smaller levels as resolution is increased. Indeed, it is well known that the resolution obtained in this particular problem is strongly related with the number of points used for the reconstructions, so that, unlike the other problems presented herein, in this case increasing arbitrarily the order of the scheme seems to improve considerably its efficiency. Finally, in order to stress out the performance of our schemes, with the different time discretizations, at a same resolution, we show in Table \ref{dmr_cpu} the computational time taken by all these combinations.

\begin{table}
  \centering
  \begin{tabular}{|c|c|c|c|}
    \hline
     & RK3 & LW & ALW \\
    \hline
    JS-WENO5 & 85.76 & 65.63 & 64.53 \\
    JS-WENO3 & 56.56 & 33.72 & 30.46 \\
    YC-WENO3 & 57.04 & 34.42 & 30.82 \\ 
    OWENO3 & 58.83  & 35.23 & 31.48 \\
    \hline
  \end{tabular}
  \caption{Example 5b: computational time (seconds) with a resolution
    of $256\times64$ grid points. $T=0.2$,
    $\textnormal{CFL}=0.25$.}
  \label{dmr_cpu}
\end{table}

One can see that, for instance, the JS-WENO5 schemes
combined with  the third-order
TVD Runge-Kutta  time discretization \cite{shuosher89} is almost three times
slower than any of the third-order optimal WENO approaches with an
approximate Lax-Wendroff time discretization.

\subsubsection*{Examples 6a and 6b: 2D Riemann problem}

Now we solve numerically a Riemann problem for the 2D Euler equations on the domain $(0,1)\times(0,1)$. An early study of Riemann problems for 2D Euler equations is \cite{SchulzRinne}. The initial data is taken as
\begin{align*}
  \boldsymbol{u}(x, y, 0)=(\rho(x, y, 0), \rho(x, y, 0)v^x(x, y, 0), \rho(x, y, 0)v^y(x, y, 0), E(x, y, 0))
  \end{align*}
with the constants (see \cite[Sect.~3, Config.~3]{KurganovTadmor}):
\[
\begin{pmatrix}
  \rho(x, y, 0)\\
  v^x(x, y, 0)\\
  v^y(x, y, 0)\\
  p(x, y, 0)
\end{pmatrix}^{\mathrm{T}} 
  =
\begin{cases}
  (1.5, 0, 0, 1.5) & \text{for $x>0.5$, $y>0.5$,}  \\
  (0.5323, 1.206, 0, 0.3) & \text{for  $x\leq0.5$, $y>0.5$,}  \\
  (0.138, 1.206, 1.206, 0.029) &  \text{for $x\leq0.5$, $y\leq0.5$,}  \\
  (0.5323, 0, 1.206, 0.3) &  \text{for $x>0.5$, $y\leq0.5$,} 
\end{cases}
\]
with the same equation of state as in the previous test.

\begin{figure}[t] 
  \centering
  \begin{tabular}{cc}
    \includegraphics[width=0.45\textwidth]{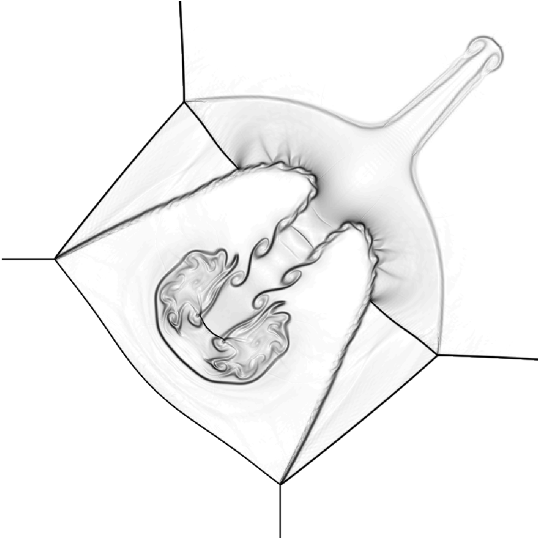} & \includegraphics[width=0.45\textwidth]{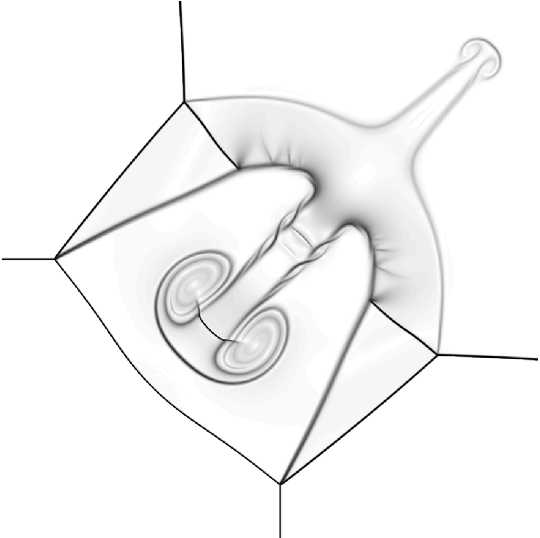} \\
    JS-WENO5 & JS-WENO3 \\
    \includegraphics[width=0.45\textwidth]{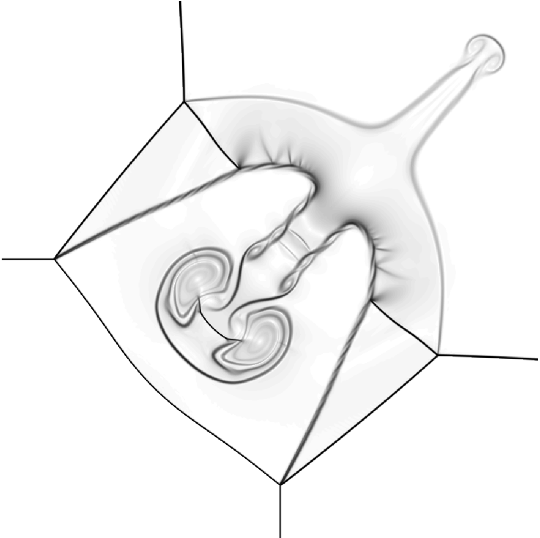} & \includegraphics[width=0.45\textwidth]{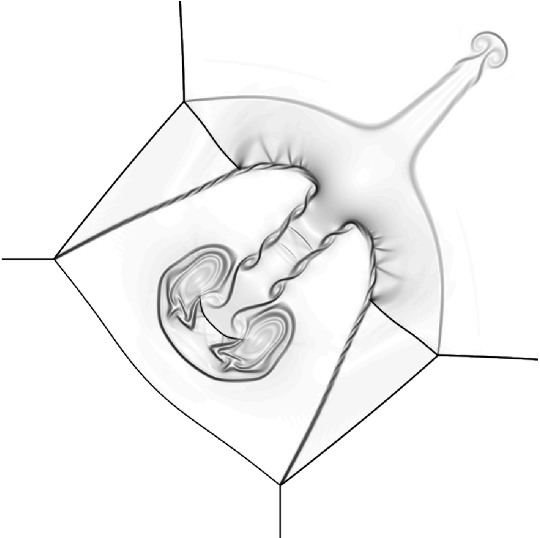} \\
    YC-WENO3 & OWENO3
  \end{tabular}
  \caption{Example 6a: 2D Riemann problem,
    $2560\times2560$. $T=0.3$. Schlieren plot of the density field.}
  \label{riemann2560}
\end{figure}

We impose outflow boundary conditions everywhere and run this test up
to time $T=0.3$. The
results can be observed in Figure \ref{riemann2560} for a
resolution of $2560\times2560$ points.


It can be seen that the order from lower to higher
resolution is again the following one: JS-WENO3, YC-WENO3, OWENO3 and
JS-WENO5, being the two latter ones close to reach other. This is very
significant
if one takes  into account that OWENO3 is faster than JS-WENO5.

\begin{figure}[t] 
  \centering
  \begin{tabular}{cc}
    \includegraphics[width=0.85\textwidth]{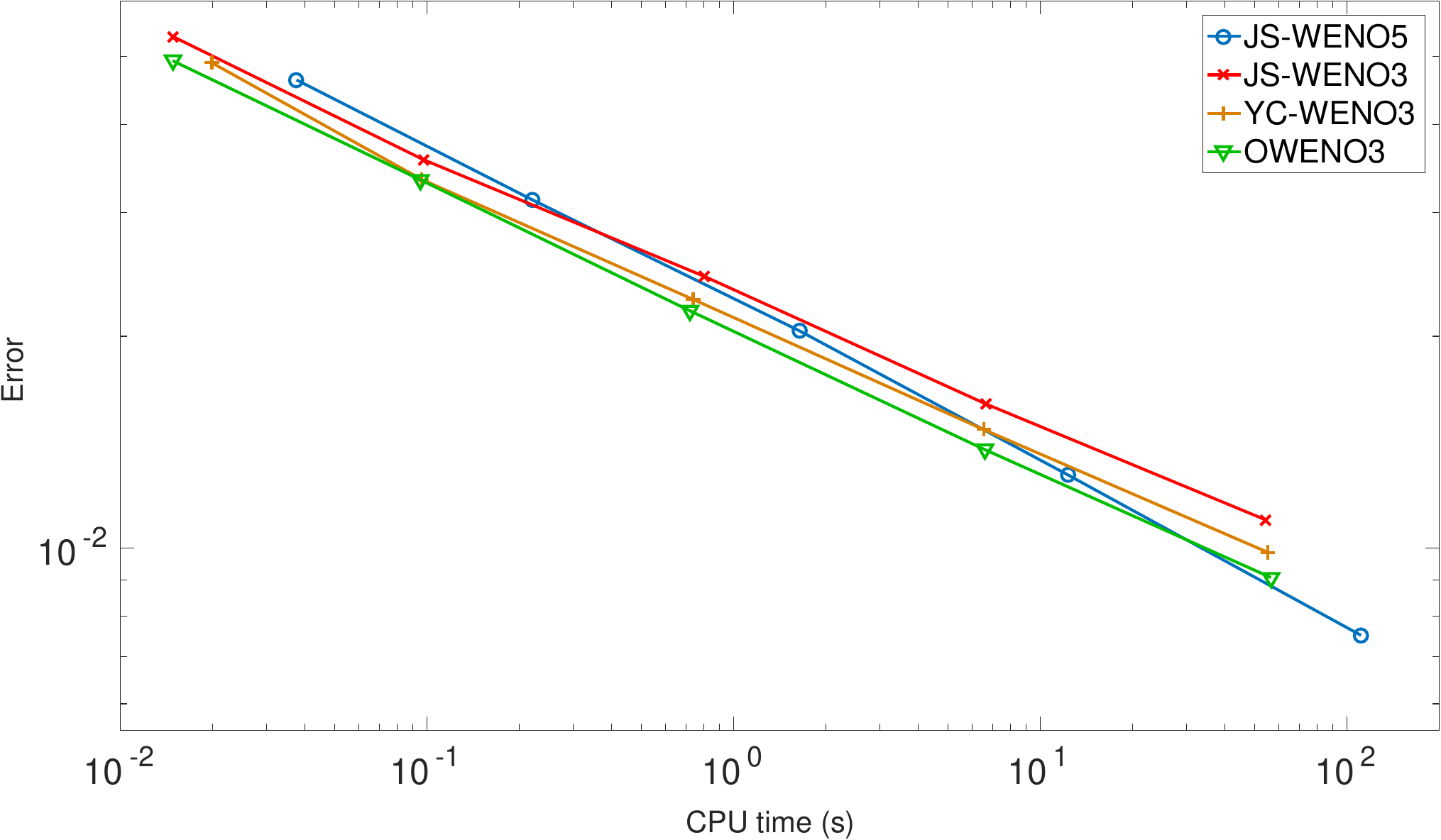}
  \end{tabular}
  \caption{Example 6b: 2D Riemann problem, efficiency plot.}
  \label{riemanneff}
\end{figure}

With the purpose of analyzing more accurately the efficiency
associated to each scheme, we now use the solutions computed with the
grid of $2560\times2560$ points as reference solutions to perform
efficiency tests by comparing error versus CPU time involving
numerical solutions with grid sizes $16\cdot2^n\times16\cdot2^n$,
$n\in\{0,1,2,3,4\}$, for the involved schemes. The results are shown in
Figure \ref{riemanneff} and again indicate a higher performance for
the OWENO3 scheme with respect to their third-order traditional
counterparts.

\subsubsection*{Examples 7a and 7b: 2D Mach 3 wind tunnel with a step}

This well-known problem involves 2D Euler equations and was proposed in
\cite{Emery68,Colella}. It consists in a wind tunnel of height 1 and width 3,
with a 0.2-height step located at 0.6 units from
the left side. A right-going Mach 3 flow is considered, such that the
initial conditions in the whole domain are $\rho=1.4$, $v^x=3$,
$v^y=0$, and $p=1$. The boundary conditions are reflecting both in the
step and the upper and bottom boundaries, inflow at the left with the
same values as the initial condition and outflow at the right.

\begin{figure}[t] 
  \centering
  \includegraphics[width=0.9\textwidth]{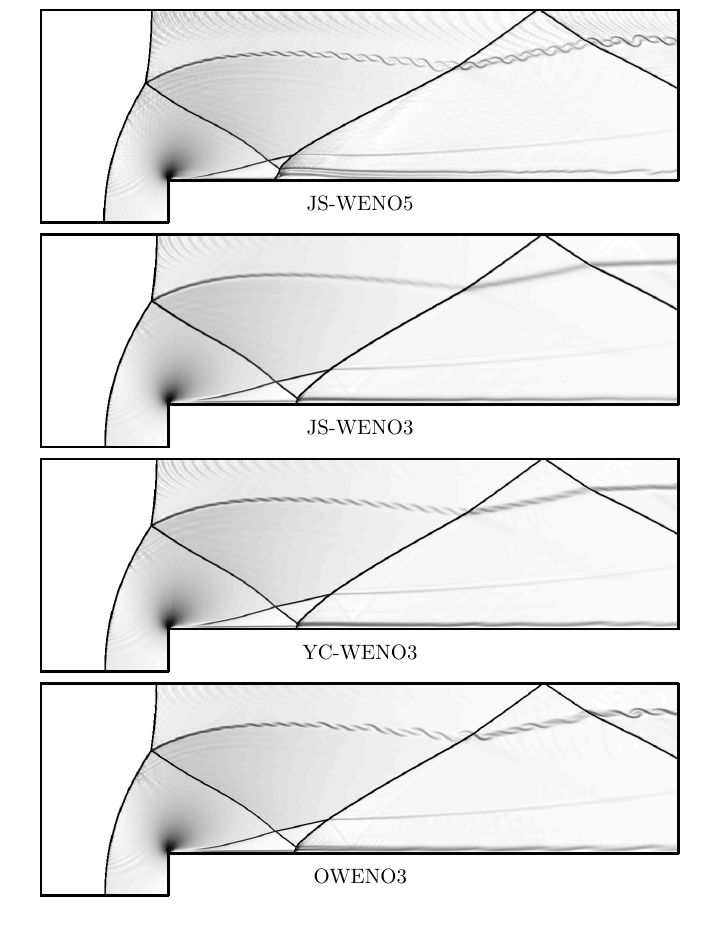}
  \caption{Example 7a: 2D Mach 3 wind tunnel with a step,
    $h_x=h_y=1/400$. $T=4$. Schlieren plot of the density field.}
  \label{step400}
\end{figure}

\begin{figure}[t] 
  \centering
  \includegraphics[width=0.9\textwidth]{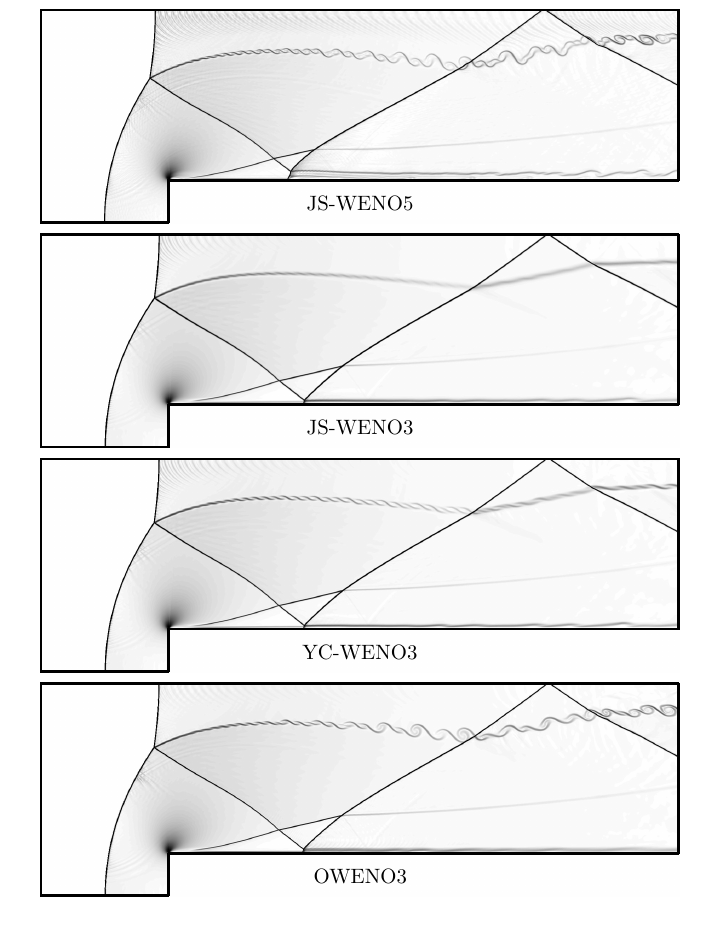}
  \caption{Example 7a: 2D Mach 3 wind tunnel with a step,
    $h_x=h_y=1/600$. $T=4$. Schlieren plot of the density field.}
  \label{step600}
\end{figure}

We perform two experiments with grid sizes $h_x=h_y=1/400$ and
$h_x=h_y=1/600$, computing the numerical
solution using the JS-WENO5, JS-WENO3, YC-WENO3 and OWENO3 schemes
until $T=4$,
which can be seen in Figures \ref{step400} and \ref{step600}. A
comparison regarding the
CPU time involving all the schemes used for a resolution of
$h_x=h_y=1/40$ points can be also found on Table \ref{step_cpu}.

\begin{table}
  \centering
  \begin{tabular}{|c|c|c|c|}
    \hline
     & RK3 & LW & ALW \\
    \hline
    JS-WENO5 & 69.46 & 55.95 & 54.60 \\
    JS-WENO3 & 52.40 & 27.29 & 26.69 \\
    YC-WENO3 & 52.71 & 27.34 & 27.00 \\ 
    OWENO3 & 53.37 & 27.65 & 27.24 \\
    \hline
  \end{tabular}
  \caption{Example 7b: computational time (seconds) with a resolution of $120\times40$ grid points. $T=4$.}
  \label{step_cpu}
\end{table}

From the results, it can be seen that OWENO3 provides a sharper
profile at the turbulent zone near the top of the domain, while
having a very similar computational cost than their classical
third-order counterparts.

\section{Conclusions}\label{sec:conclusions}

In this paper it has been proven that a third-order
    interpolator with a 3-points stencil cannot simultaneously detect
    discontinuities and keep the optimal third-order accuracy near
    critical points unless a scale-dependent parameter is
    used. As a
    consequence, a third-order
scheme, whose numerical flux interpolator includes a
    fourth additional node (used only for the computation of the
    weights),  based on a
WENO approach with unconditionally third-order optimal
accuracy on smooth data, and without relying on any tuning parameter,
has been presented. The resulting scheme maintains the
    width of the domain of dependence (a
stencil of at most 4 points is used to obtain each numerical flux, as in
the traditional third-order WENO schemes) and the accuracy properties
of the proposed method have been proved theoretically and confirmed numerically along
experiments involving algebraic problems and hyperbolic
conservation laws. The novel scheme is more efficient than the other 
three-order methods considered and in most cases outperforms
 even the classical 
fifth-order JS-WENO scheme. Only in some problems
involving very small-scale features, like the double Mach reflection
test, the fifth-order method is competitive. However, it must be also taken into account that the third-order
schemes considered in this paper allow higher values of the CFL in the
aforementioned complicated problems, and therefore, even in these cases third-order schemes may be worth
being used instead as well.

\section*{Acknowledgments} 

AB, PM and DZ are supported by Spanish MINECO pro\-ject
MTM2017-83942-P.  
RB is supported by 
 CONICYT/PIA/AFB170001; CRHIAM, pro\-ject CONICYT/FON\-DAP/15130015; Fondecyt
pro\-ject 1170473; and the \mbox{INRIA} Associated Team ``Efficient numerical schemes for non-local transport phenomena'' (NOLOCO; 2018--2020). 
PM is also
supported by Conicyt (Chile),   pro\-ject PAI-MEC, folio 80150006.  
 DZ is also supported by Conicyt (Chile) through Fondecyt project 3170077.

\end{document}